# CONTINUOUS-TIME GARCH PROCESSES

BY PETER BROCKWELL,[1] ERDENEBAATAR CHADRAA
AND ALEXANDER LINDNER[2]

*Colorado State University, Colorado State University and Technical University of Munich*

A family of continuous-time generalized autoregressive conditionally heteroscedastic processes, generalizing the COGARCH$(1,1)$ process of Klüppelberg, Lindner and Maller [*J. Appl. Probab.* **41** (2004) 601–622], is introduced and studied. The resulting COGARCH$(p,q)$ processes, $q \geq p \geq 1$, exhibit many of the characteristic features of observed financial time series, while their corresponding volatility and squared increment processes display a broader range of autocorrelation structures than those of the COGARCH$(1,1)$ process. We establish sufficient conditions for the existence of a strictly stationary nonnegative solution of the equations for the volatility process and, under conditions which ensure the finiteness of the required moments, determine the autocorrelation functions of both the volatility and the squared increment processes. The volatility process is found to have the autocorrelation function of a continuous-time autoregressive moving average process.

**1. Introduction.** In financial econometrics, discrete-time GARCH (generalized autoregressive conditionally heteroscedastic) processes are widely used to model the returns at regular intervals on stocks, currency investments and other assets. Specifically, a GARCH process $(\xi_n)_{n \in \mathbb{N}}$ typically represents the increments $\ln P_n - \ln P_{n-1}$ of the logarithms of the asset price at times $1, 2, 3, \ldots$. These models capture many of the so-called *stylized features* of such data, for example, tail heaviness, volatility clustering and dependence without correlation. For GARCH processes with finite fourth

Received April 2005; revised September 2005.
[1] Supported by NSF Grant DMS-03-08109.
[2] Supported by Deutsche Forschungsgemeinschaft Grant Li 1026/2-1.

moments, the autocorrelation functions of both the squared process and the associated volatility process are those of autoregressive moving average (ARMA) processes. The squared GARCH(1,1) process, for example, has the autocorrelation function of an ARMA(1,1) process and the corresponding volatility has the autocorrelation function of an AR(1) process. Higher-order GARCH$(p,q)$ processes were introduced to allow for the possibility of a broader range of autocorrelations for the volatility and the squared increment processes.

Various attempts have been made to capture the stylized features of financial time series using continuous-time models. The interest in continuous-time models stems from their use in modeling irregularly spaced data, their use in financial applications such as option pricing and the current widespread availability of high-frequency data. In continuous time it is natural to model the logarithm of the asset price itself, that is, $G_t = \ln P_t$, rather than its increments as in discrete time.

Notable among these attempts is the GARCH diffusion approximation of Nelson [23]. (See also [12] and [13].) Although the GARCH process is driven by a single noise sequence, the diffusion limit is driven by two independent Brownian motions $(W_t^{(1)})_{t \geq 0}$ and $(W_t^{(2)})_{t \geq 0}$. For example, the GARCH(1,1) diffusion limit satisfies

$$
\begin{aligned}
dG_t &= \sigma_t \, dW_t^{(1)}, \\
d\sigma_t^2 &= \theta(\gamma - \sigma_t^2) + \rho \sigma_t^2 \, dW_t^{(2)}, \qquad t \geq 0.
\end{aligned}
\tag{1.1}
$$

The behavior of this diffusion limit is therefore rather different from that of the GARCH process itself since the volatility process $(\sigma_t^2)_{t \geq 0}$ evolves independently of the driving process $(W_t^{(1)})_{t \geq 0}$ in the first of the equations (1.1).

Another approach is via the stochastic volatility model of Barndorff-Nielsen and Shephard [3, 4] in which the volatility process $\sigma^2$ is an Ornstein–Uhlenbeck (O–U) process driven by a nondecreasing Lévy process and $G$ satisfies an equation of the form $dG_t = \mu \, dt + \sigma_t \, dW_t$, where $W$ is a Brownian motion independent of the Lévy process. The autocorrelation function of the Lévy-driven O–U volatility process has the form $\rho(h) = \exp(-c|h|)$ for some $c > 0$, but this class can be extended by specifying the volatility to be a superposition of O–U processes as in [2] or a Lévy-driven CARMA (continuous-time ARMA) process as in [10]. As in Nelson's diffusion limit, the process $G$ is again driven by two independent noise processes and the volatility process $\sigma^2$ evolves independently of the process $W$ in the equation for $G$.

A continuous-time analog of the GARCH(1,1) process, denoted COGARCH(1,1), was recently constructed and studied by Klüppelberg,

Lindner and Maller [19]. Their construction is based on taking a limit of an explicit representation of the discrete-time GARCH(1, 1) process to obtain a continuous-time analog. Since no such representation exists for higher-order discrete-time GARCH processes, a different approach is needed to construct higher-order continuous-time analogs. In this paper we do this by specifying a system of Lévy-driven stochastic differential equations for the processes $G$ and $\sigma^2$. If the volatility process $\sigma^2$ is strictly stationary we refer to the process $G$ as a COGARCH$(p, q)$ process. In the special case $p = q = 1$ we recover the COGARCH(1, 1) process of Klüppelberg, Lindner and Maller [19]. In general we obtain a class of processes $G$ with uncorrelated increments but for which the corresponding volatility and squared increment processes exhibit a broad range of autocorrelation functions. The volatility process has the autocorrelation function of a continuous-time ARMA process.

The construction of the COGARCH(1, 1) process due to Klüppelberg, Lindner and Maller [19] starts from the defining equations of the discrete-time GARCH(1, 1) process $(\xi_n)_{n \in \mathbb{N}_0}$,

(1.2)
$$\xi_n = \varepsilon_n \sigma_n,$$
$$\sigma_n^2 = \alpha_0 + \alpha_1 \xi_{n-1}^2 + \beta_1 \sigma_{n-1}^2, \qquad n \in \mathbb{N}_0,$$

where $\alpha_0, \alpha_1$ and $\beta_1$ are all strictly positive, and $(\varepsilon_n)_{n \in \mathbb{N}_0}$ is a sequence of i.i.d. (independent and identically distributed) random variables with mean 0 and variance 1. The recursions (1.2) can be solved to give

$$\sigma_n^2 = \left(\sigma_0^2 + \alpha_0 \int_0^n \exp\left[-\sum_{j=0}^{\lfloor s \rfloor} \log(\beta_1 + \alpha_1 \varepsilon_j^2)\right] ds\right)$$
$$\times \exp\left[\left(\sum_{j=0}^{n-1} \log(\beta_1 + \alpha_1 \varepsilon_j^2)\right)\right],$$

where $\lfloor s \rfloor$ denotes the integer part of $s \in \mathbb{R}$. The COGARCH(1, 1) equations are then obtained by replacing the driving noise sequence $(\varepsilon_n)_{n \in \mathbb{N}_0}$ by the jumps $(\Delta L_t = L_t - L_{t-})_{t \geq 0}$ of a Lévy process. More precisely, observing that $\sum_{j=0}^{n-1} \log(\beta_1 + \alpha_1 \varepsilon_j^2) = n \log \beta_1 + \sum_{j=0}^{n-1} \log(1 + (\alpha_1/\beta_1)\varepsilon_j^2)$ for $\beta_1 > 0$ and writing $\eta$ for $-\log \beta_1$, $\omega_0$ for $\alpha_0$ and $\omega_1$ for $\alpha_1$, leads to the equations

(1.3) $\quad dG_t = \sigma_t\, dL_t, \qquad\qquad\qquad\qquad t > 0, \qquad G_0 = 0,$

(1.4) $\quad \sigma_t^2 = \left(\sigma_0^2 + \omega_0 \int_0^t e^{X_s}\, ds\right) e^{-X_{t-}}, \qquad t \geq 0,$

where

(1.5) $\quad X_t := \eta t - \sum_{0 < s \leq t} \log(1 + \omega_1 e^{\eta}(\Delta L_s)^2).$

Here, $\omega_0 > 0, \omega_1 \geq 0, \eta > 0$ and $\sigma_0^2$ is independent of $(L_t)_{t \geq 0}$. The COGARCH$(1,1)$ process is the solution $G$ of these equations and, under specified conditions on the coefficients and the distribution of $\sigma_0^2$, the volatility process $\sigma^2$ is strictly stationary and $G$ has stationary increments.

The COGARCH$(1,1)$ process with stationary volatility has been shown to have many of the features of the discrete time GARCH$(1,1)$ process. As shown in [19, 20], the COGARCH$(1,1)$ process has uncorrelated increments, while the autocorrelation functions of the volatility $\sigma^2$ and of the squared increments of $G$ decay exponentially. Furthermore, the COGARCH$(1,1)$ process has heavy tails and volatility clusters at high levels; see [14] and [20]. For an overview of the extremes of stochastic volatility models, see [14] and [21]. Also, observe that many of the features of the COGARCH$(1,1)$ process can be obtained in a more general setting, as in [21].

In the next section we specify a system of stochastic differential equations for the COGARCH$(p,q)$ process $G$ and its associated volatility process, which we shall denote by $V$. This is directly motivated by the corresponding structure of the discrete-time GARCH$(p,q)$ process. We then show that the solution of these equations coincides with that of the COGARCH$(1,1)$ equations if $p = q = 1$. Notation and definitions used throughout the paper are given at the end of Section 2.

In Section 3 we give sufficient conditions for the existence of a strictly stationary volatility process. In the COGARCH$(1,1)$ case, these are exactly the necessary and sufficient conditions obtained by Klüppelberg, Lindner and Maller [19, 20]. More detailed results are given in the special case when the driving Lévy process is compound Poisson. The proofs rely on the fact that the state vector of the COGARCH$(p,q)$ process, sampled at uniformly spaced discrete times, satisfies a multivariate random recurrence equation.

In Section 4 we focus on the autocorrelation structure of the stationary volatility process. Just as the discrete-time GARCH volatility process has the autocorrelation function of an ARMA process, the COGARCH volatility process has the autocorrelation function of a CARMA process.

Section 5 deals with conditions which ensure positivity of the volatility, while the autocorrelation structure of the squared increments of the COGARCH process itself is obtained in Section 6. The results are illustrated with simulations in Section 7. So as not to disturb the flow of the arguments, proofs of the results are postponed to Sections 8–11.

**2. The COGARCH$(p,q)$ equations.** Let $(\varepsilon_n)_{n \in \mathbb{N}_0}$ be an i.i.d. sequence of random variables. Let $p, q \geq 0$. Then the GARCH$(p,q)$ process $(\xi_n)_{n \in \mathbb{N}_0}$ is defined by the equations

$$\xi_n = \sigma_n \varepsilon_n,$$
(2.1) $$\sigma_n^2 = \alpha_0 + \alpha_1 \xi_{n-1}^2 + \cdots + \alpha_p \xi_{n-p}^2$$

$$+ \beta_1 \sigma_{n-1}^2 + \cdots + \beta_q \sigma_{n-q}^2, \qquad n \geq s,$$

where $s := \max(p,q)$, $\sigma_0^2, \ldots, \sigma_{s-1}^2$ are i.i.d. and independent of the i.i.d. sequence $(\varepsilon_n)_{n \geq s}$, and $\xi_n = G_{n+1} - G_n$ represents the increment at time $n$ of the log asset price process $(G_n)_{n \in \mathbb{N}_0}$. Note that the continuous-time GARCH process will be a model for $(G_t)_{t \geq 0}$ and not for its increments as in discrete time.

Equation (2.1) shows that the volatility process $(\sigma_n^2)_{n \in \mathbb{N}_0}$ can be viewed as a "self-exciting" $\mathrm{ARMA}(q, p-1)$ process driven by the noise sequence $(\sigma_{n-1}^2 \varepsilon_{n-1}^2)_{n \in \mathbb{N}}$. Motivated by this observation, we will define a continuous-time GARCH model for the log asset price process $(G_t)_{t \geq 0}$ of order $(p,q)$ by

$$dG_t = \sigma_t \, dL_t, \qquad t > 0, \qquad G_0 = 0,$$

where $(\sigma_t^2)_{t \geq 0}$ is a $\mathrm{CARMA}(q, p-1)$ process driven by a suitable replacement for the discrete time driving noise sequence $(\sigma_{n-1}^2 \varepsilon_{n-1}^2)_{n \in \mathbb{N}}$.

The state–space representation of a Lévy-driven $\mathrm{CARMA}(q, p-1)$ process $(\psi_t)_{t \geq 0}$ with driving Lévy process $L$, location parameter $c$, moving average coefficients $\alpha_1, \ldots, \alpha_p$, autoregressive coefficients $\beta_1, \ldots, \beta_q$ and $q \geq p$ is (see [8])

$$\psi_t = c + \mathbf{a}' \zeta_t,$$

$$d\zeta_t = \begin{bmatrix} 0 & 1 & 0 & \cdots & 0 \\ 0 & 0 & 1 & \cdots & 0 \\ \vdots & \vdots & \vdots & \ddots & 0 \\ 0 & 0 & 0 & \cdots & 1 \\ -\beta_q & -\beta_{q-1} & -\beta_{q-2} & \cdots & -\beta_1 \end{bmatrix} \zeta_t \, dt + \begin{bmatrix} 0 \\ 0 \\ \vdots \\ 0 \\ 1 \end{bmatrix} dL_t,$$

where $\mathbf{a}' = [\alpha_1, \ldots, \alpha_q]$, $\alpha_j := 0$ for $j > p$ and the coefficient matrix in the last equation is $-\beta_1$ if $q = 1$. (The $\mathrm{CARMA}(q, p-1)$ process $(\psi_t)_{t \geq 0}$ is a strictly stationary solution of these equations, which exists under conditions found in [9].) To obtain a continuous-time analog of the equation (2.1), we suppose that the volatility process $(\sigma_t^2)_{t \geq 0}$ has the state–space representation of a $\mathrm{CARMA}(q, p-1)$ process in which the driving Lévy process $(L_t)$ is replaced by a continuous-time analog of the driving process $(\sigma_{n-1}^2 \varepsilon_{n-1}^2)_{n \in \mathbb{N}}$ in (2.1).

The increments of the driving process in continuous time should correspond to the increments of the discrete-time process:

$$R_n^{(d)} := \sum_{i=0}^{n-1} \xi_i^2 = \sum_{i=0}^{n-1} \sigma_i^2 \varepsilon_i^2.$$

We therefore replace the innovations $\varepsilon_n$ by the jumps $\Delta L_t$ of a Lévy process $(L_t)_{t \geq 0}$ to obtain the continuous-time analog

$$R_t := \sum_{0 < s \leq t} \sigma_{s-}^2 (\Delta L_s)^2, \qquad t > 0.$$

If $L$ has no Gaussian part [i.e., $\tau_L^2 = 0$ in (2.2) below], we recognize $R$ as the quadratic covariation of $G$, that is,

$$R_t = \sum_{0<s\leq t} \sigma_{s-}^2 (\Delta L_s)^2 = \int_0^t \sigma_{s-}^2 \, d[L,L]_s = [G,G]_t.$$

If $L$ has a Gaussian part, then still $\sum_{0<s\leq t}(\Delta L_s)^2 = [L,L]^{(d)}$, the *discrete part of the quadratic covariation*, and we have in general

$$R_t = \int_0^t \sigma_{s-}^2 \, d[L,L]_s^{(d)}, \quad \text{that is,} \quad dR_t = \sigma_{t-}^2 \, d[L,L]_t^{(d)}.$$

Recall that for a Lévy process $L = (L_t)_{t\geq 0}$ the characteristic function $E(e^{i\theta L_t}), \theta \in \mathbb{R}$, can be written in the form

$$\begin{aligned}
E(e^{i\theta L_t}) = \exp\bigg( t\bigg( i\gamma_L \theta - \tau_L^2 \frac{\theta^2}{2} \\
+ \int_\mathbb{R} (e^{i\theta x} - 1 - i\theta x \mathbf{1}_{|x|\leq 1}) \, d\nu_L(x) \bigg) \bigg).
\end{aligned} \tag{2.2}$$

The constants $\gamma_L \in \mathbb{R}$, $\tau_L^2 \geq 0$ and the measure $\nu_L$ on $\mathbb{R}$ form the *characteristic triplet* of $L$. As usual, the Lévy measure $\nu_L$ is required to satisfy $\int_\mathbb{R} \min(1, x^2) \, d\nu_L(x) < \infty$. For more information on Lévy processes, we refer to the books by Applebaum [1], Bertoin [5] or Sato [25].

The COGARCH$(p,q)$ equations will now be obtained by specifying that the volatility process $V\ (= \sigma^2)$ should satisfy continuous-time ARMA equations driven by the process $R$ defined above. Provided $V$ is nonnegative almost surely (conditions for which are given in Section 5), we can define a process $G$ by the equations $G_0 = 0$ and $dG_t = \sqrt{V_t} \, dL_t$. Under conditions ensuring that $V$ is also strictly stationary, we refer to $G$ as a COGARCH$(p,q)$ process. As we shall see, when $p = q = 1$, the solution of the COGARCH equations coincides with that of the COGARCH$(1,1)$ equations (1.3)–(1.5) of [19]. [The parameters $\beta_1, \ldots, \beta_q$ and $\alpha_1, \ldots, \alpha_p$ in the following definition should not be confused with the parameters denoted by the same symbols in the defining equation (2.1) of the discrete-time GARCH process.]

DEFINITION 2.1 [The COGARCH$(p,q)$ equations]. If $p$ and $q$ are integers such that $q \geq p \geq 1$, $\alpha_0 > 0$, $\alpha_1, \ldots, \alpha_p \in \mathbb{R}$, $\beta_1, \ldots, \beta_q \in \mathbb{R}$, $\alpha_p \neq 0$, $\beta_q \neq 0$ and $\alpha_{p+1} = \cdots = \alpha_q = 0$, we define the $(q \times q)$ matrix $B$ and the vectors $\mathbf{a}$ and $\mathbf{e}$ by

$$B = \begin{bmatrix} 0 & 1 & 0 & \cdots & 0 \\ 0 & 0 & 1 & \cdots & 0 \\ \vdots & \vdots & \vdots & \ddots & \vdots \\ 0 & 0 & 0 & \cdots & 1 \\ -\beta_q & -\beta_{q-1} & -\beta_{q-2} & \cdots & -\beta_1 \end{bmatrix},$$

$$\mathbf{a} = \begin{bmatrix} \alpha_1 \\ \alpha_2 \\ \vdots \\ \alpha_{q-1} \\ \alpha_q \end{bmatrix}, \qquad \mathbf{e} = \begin{bmatrix} 0 \\ 0 \\ \vdots \\ 0 \\ 1 \end{bmatrix},$$

with $B := -\beta_1$ if $q = 1$. Then if $L = (L_t)_{t \geq 0}$ is a Lévy process with nontrivial Lévy measure, we define the (left-continuous) *volatility process* $V = (V_t)_{t \geq 0}$ with parameters $B$, $\mathbf{a}$, $\alpha_0$ and driving Lévy process $L$ by

$$V_t = \alpha_0 + \mathbf{a}' \mathbf{Y}_{t-}, \qquad t > 0, \qquad V_0 = \alpha_0 + \mathbf{a}' \mathbf{Y}_0,$$

where the *state process* $\mathbf{Y} = (\mathbf{Y}_t)_{t \geq 0}$ is the unique cadlag solution of the stochastic differential equation

(2.3) $$d\mathbf{Y}_t = B \mathbf{Y}_{t-} \, dt + \mathbf{e}(\alpha_0 + \mathbf{a}' \mathbf{Y}_{t-}) \, d[L,L]_t^{(d)}, \qquad t > 0,$$

with initial value $\mathbf{Y}_0$, independent of the driving Lévy process $(L_t)_{t \geq 0}$. If the process $(V_t)_{t \geq 0}$ is strictly stationary and nonnegative almost surely, we say that $G = (G_t)_{t \geq 0}$, given by

$$dG_t = \sqrt{V_t} \, dL_t, \qquad t > 0, \qquad G_0 = 0,$$

is a COGARCH$(p,q)$ *process* with parameters $B$, $\mathbf{a}$, $\alpha_0$ and driving Lévy process $L$.

That there is in fact a unique solution of (2.3) for any starting random vector $\mathbf{Y}_0$ follows from standard theorems on stochastic differential equations (e.g., [24], Chapter V, Theorem 7). The stochastic integrals are interpreted with respect to the filtration $\mathbb{F} = (\mathcal{F}_t)_{t \geq 0}$, which is defined to be the smallest right-continuous filtration such that $\mathcal{F}_0$ contains all the $P$-null sets of $\mathcal{F}$, $(L_t)_{t \geq 0}$ is adapted and $\mathbf{Y}_0$ is $\mathcal{F}_0$-measurable.

Without restrictions on $\alpha_0, \mathbf{a}$ and $B$, the process $V$ is not necessarily nonnegative, in which case $G$ is undefined. Conditions which ensure that $V$ is nonnegative will be discussed in Section 5. In particular, it will be shown that if $\mathbf{a}' e^{Bt} \mathbf{e} \geq 0$ for all $t \geq 0$ and $\mathbf{Y}_0$ is such that $V$ is strictly stationary, then $V$ is nonnegative with probability 1. Even if $V$ takes negative values, however, the process is of some interest in its own right and many of our results for $V$ are valid without the nonnegativity restriction.

Conditions for stationarity of $V$ are discussed in Section 3.

We next show that if $p = q = 1$, the solution of the COGARCH equations in Definition 2.1 coincides with the solution of the COGARCH$(1,1)$ equations of Klüppelberg, Lindner and Maller [19].

THEOREM 2.2. *Suppose that $p = q = 1$, and that $\alpha_0, \alpha_1$ and $\beta$ are all strictly positive. Then the processes $(G_t)_{t \geq 0}$ and $(V_t)_{t \geq 0}$ of Definition 2.1 are, respectively, the processes $(G_t)_{t \geq 0}$ and $(\sigma_t^2)_{t \geq 0}$ defined by (1.3)–(1.5), with parameters $\omega_0 = \alpha_0 \beta_1$, $\omega_1 = \alpha_1 e^{-\beta_1}$ and $\eta = \beta_1$.*

Proof. From
$$d\mathbf{Y}_t = -\beta_1 \mathbf{Y}_t\, dt + V_t\, d[L,L]_t^{(d)} \quad \text{and} \quad V_{t+} = \alpha_0 + \alpha_1 \mathbf{Y}_t$$
it follows that
$$dV_{t+} = \alpha_1\, d\mathbf{Y}_t = -\alpha_1 \beta_1 \frac{V_t - \alpha_0}{\alpha_1}\, dt + \alpha_1 V_t\, d[L,L]_t^{(d)},$$
and hence that
$$V_{t+} = \alpha_0 \beta_1 t - \beta_1 \int_0^t V_s\, ds + \alpha_1 \sum_{0 < s \leq t} V_s (\Delta L_s)^2 + V_0.$$
However, this equation is also satisfied by the volatility process $(\sigma_t^2)_{t \geq 0}$ of (1.4) when $\omega_0 = \alpha_0 \beta_1$, $\eta = \beta_1$ and $\omega_1 = \alpha_1 e^{-\beta_1}$, as shown in Proposition 3.2 of [19], and uniqueness of the solution gives the claim. □

We conclude this section with a few definitions and some notation which will be used throughout the paper.

DEFINITION 2.3. Let $\mathbf{a}$ and $B$ be as in Definition 2.1. Then the *characteristic polynomials* associated with $\mathbf{a}$ and $B$ are given by
$$a(z) := \alpha_1 + \alpha_2 z + \cdots + \alpha_p z^{p-1}, \qquad z \in \mathbb{C},$$
$$b(z) := z^q + \beta_1 z^{q-1} + \cdots + \beta_q, \qquad z \in \mathbb{C}.$$
The eigenvalues of the matrix $B$ (which are exactly the zeroes of $b$) will be denoted by $\lambda_1, \ldots, \lambda_q$ and assumed to be ordered in such a way that
$$\Re \lambda_q \leq \Re \lambda_{q-1} \leq \cdots \leq \Re \lambda_1$$
(where $\Re \lambda_i$ denotes the real part of $\lambda_i$). Furthermore, define
$$\lambda := \lambda(B) := \Re \lambda_1.$$

For the rest of the paper, convergence in probability will be denoted by "$\xrightarrow{P}$", uniform convergence on compacts in probability by "$\xrightarrow{ucp}$" and equality in distribution by "$\stackrel{d}{=}$". For $x \in \mathbb{R}$ we shall write $\log^+(x)$ for $\log(\max\{1,x\})$. The transpose of a column vector $\mathbf{c} \in \mathbb{C}^q$ will be denoted by $\mathbf{c}'$. If $\|\cdot\|$ is a vector norm in $\mathbb{C}^q$, then the natural matrix norm of the $(q \times q)$ matrix $C$ is defined as $\|C\| = \sup_{\mathbf{c} \in \mathbb{C}^q \setminus \{0\}} \frac{\|C\mathbf{c}\|}{\|\mathbf{c}\|}$. Correspondingly, for $r \in [1, \infty]$, we denote by $\|\cdot\|_r$ both the vector $L^r$-norm and the associated natural matrix norm. Recall that the natural matrix norms of the $L^1, L^2$ and $L^\infty$ vector norms are the column-sum norm, the spectral norm and the row-sum norm, respectively.

The $(q \times q)$ identity matrix will be denoted by $I_q$ or simply $I$, and the canonical vector $(0, \ldots, 0, 1, 0, \ldots, 0)'$, with $i$th component equal to 1, will be denoted by $\mathbf{e}_i$. For $\mathbf{e}_q$ we simply write $\mathbf{e}$. By $\operatorname{diag}(\lambda_1, \ldots, \lambda_q)$ we mean the diagonal $(q \times q)$ matrix with these entries on the diagonal. The Kronecker product of two $(q \times q)$ matrices $A$ and $B$ will be denoted by $A \otimes B$, and by $\operatorname{vec}(A)$ we denote the column vector in $\mathbb{C}^{q^2}$ which arises from $A$ by stacking the columns of $A$ in a vector (starting with the first column). For the properties of the Kronecker product, we refer to Lütkepohl [22].

**3. Stationarity conditions.** In this section we consider conditions under which the volatility process $(V_t)_{t \geq 0}$ specified in Definition 2.1 is strictly stationary. The parameters $B$, $\mathbf{a}$ and $\alpha_0$, and the state process $(\mathbf{Y}_t)_{t \geq 0}$ are as specified in Definition 2.1. The condition (3.2) established in Theorem 3.1 below is necessary and sufficient for stationarity in the special case $p = q = 1$. For larger values of $p$ and $q$ it is sufficient only, but not unduly restrictive, since there is a rich class of models satisfying the condition. Without serious loss of generality we shall assume that the matrix $B$ can be diagonalized. Since the only eigenvectors corresponding to the eigenvalue $\lambda_i$ are constant multiples of $[1, \lambda_i, \lambda_i^2, \ldots, \lambda_i^{q-1}]'$, this is equivalent to the assumption that the eigenvalues of $B$ are distinct. Let $S$ be a matrix such that $S^{-1}BS$ is a diagonal matrix, for example,

$$(3.1) \qquad S = \begin{bmatrix} 1 & \cdots & 1 \\ \lambda_1 & \cdots & \lambda_q \\ \vdots & \cdots & \vdots \\ \lambda_1^{q-1} & \cdots & \lambda_q^{q-1} \end{bmatrix}.$$

[For this particular choice, $S^{-1}BS = \operatorname{diag}(\lambda_1, \ldots, \lambda_q)$.]

THEOREM 3.1. *Let $(\mathbf{Y}_t)_{t \geq 0}$ be the state process of the $\operatorname{COGARCH}(p, q)$ process with parameters $B$, $\mathbf{a}$ and $\alpha_0$. Suppose that all the eigenvalues of $B$ are distinct. Let $L$ be a Lévy process with nontrivial Lévy measure $\nu_L$ and suppose there is some $r \in [1, \infty]$ such that*

$$(3.2) \qquad \int_{\mathbb{R}} \log(1 + \|S^{-1}\mathbf{e}\mathbf{a}'S\|_r y^2)\, d\nu_L(y) < -\lambda = -\lambda(B)$$

*for some matrix $S$ such that $S^{-1}BS$ is diagonal. Then $\mathbf{Y}_t$ converges in distribution to a finite random variable $\mathbf{Y}_\infty$, as $t \to \infty$. It follows that if $\mathbf{Y}_0 \stackrel{d}{=} \mathbf{Y}_\infty$, then $(\mathbf{Y}_t)_{t \geq 0}$ and $(V_t)_{t \geq 0}$ are strictly stationary.*

REMARK 3.2. (a) If $(V_t)_{t \geq 0}$ is the volatility of a $\operatorname{COGARCH}(1, 1)$ process with parameters $B = -\beta_1 < 0$, $\alpha_0 > 0$ and $\alpha_1 > 0$, then $\|S^{-1}\mathbf{e}\mathbf{a}'S\|_r = \alpha_1$ and, as already indicated, the condition (3.2) is necessary and sufficient for the existence of a strictly stationary COGARCH$(1,1)$ volatility process. (See [19], Theorem 3.1.)

(b) For general $q \geq 2$, the quantity $\|S^{-1}\mathbf{ea}S\|_r$ depends on the specific choice of $S$ and on $r$. Observe that it is sufficient to find *some* $S$ and *some* $r$ such that (3.2) holds.

The proof of Theorem 3.1 will make heavy use of the general theory of multivariate random recurrence equations, as discussed by Bougerol and Picard [7], Kesten [18] and Brandt [6] (in the one-dimensional case). The COGARCH state vector satisfies such a multivariate random recurrence equation, as indicated in the following theorem.

THEOREM 3.3. *Let* $(\mathbf{Y}_t)_{t \geq 0}$ *be the state process of the* COGARCH$(p,q)$ *process with parameters* $B$, $\mathbf{a}$ *and* $\alpha_0$, *and driving Lévy process* $L$. *Then there exists a family* $(J_{s,t}, \mathbf{K}_{s,t})_{0 \leq s \leq t}$ *of random* $(q \times q)$ *matrices* $J_{s,t}$ *and random vectors* $\mathbf{K}_{s,t}$ *in* $\mathbb{R}^q$ *such that*

$$\mathbf{Y}_t = J_{s,t}\mathbf{Y}_s + \mathbf{K}_{s,t}, \qquad 0 \leq s \leq t. \tag{3.3}$$

*Furthermore, the distribution of* $(J_{s,t}, \mathbf{K}_{s,t})$ *depends only on* $t-s$, $(J_{s_1,t_1}, \mathbf{K}_{s_1,t_1})$ *and* $(J_{s_2,t_2}, \mathbf{K}_{s_2,t_2})$ *are independent for* $0 \leq s_1 \leq t_1 \leq s_2 \leq t_2$; *and for* $0 \leq s \leq u \leq t$,

$$J_{s,t} = J_{u,t} J_{s,u}. \tag{3.4}$$

*If additionally the conditions of Theorem* 3.1 *hold, then the distribution of the vector* $\mathbf{Y}_\infty$ *is, for any* $h > 0$, *the unique solution of the random fixed point equation*

$$\mathbf{Y}_\infty \stackrel{d}{=} J_{0,h}\mathbf{Y}_\infty + \mathbf{K}_{0,h}, \tag{3.5}$$

*with* $\mathbf{Y}_\infty$ *independent of* $(J_{0,h}, \mathbf{K}_{0,h})$ *on the right-hand side of* (3.5).

REMARK 3.4. (a) The stationarity condition (3.2) is easy to check. However, as the proofs of Theorems 3.1 and 3.3 show, a weaker stationary condition is the existence of a vector norm $\|\cdot\|$ and $t_0 > 0$ such that $J_{0,t_0}$ and $\mathbf{K}_{0,t_0}$ satisfy the conditions

$$E \log \|J_{0,t_0}\| < 0 \quad \text{and} \quad E \log^+ \|\mathbf{K}_{0,t_0}\| < \infty. \tag{3.6}$$

By (3.4), $E \log \|J_{0,t_0}\| < 0$ is equivalent to the requirement that there is a strictly positive value of $t_1$ such that the Lyapunov exponent of the i.i.d. sequence $(J_{t_1 n, t_1(n+1)})_{n \in \mathbb{N}_0}$, that is,

$$\lim_{n \to \infty} \frac{1}{n} E(\log \|J_{t_1(n-1), t_1 n} \cdots J_{0, t_1}\|)$$
$$= \inf_{n \in \mathbb{N}} \left( \frac{1}{n} E(\log \|J_{t_1(n-1), t_1 n} \cdots J_{0, t_1}\|) \right),$$

(which is independent of the specific norm), is strictly negative. As shown by Bougerol and Picard [7], provided $E\log^+ \|J_{0,t_1}\| < \infty$, $E\log^+ \|\mathbf{K}_{0,t_1}\| < \infty$ and a certain irreducibility condition holds, then strict negativity of the Lyapunov exponent is not only sufficient, but also necessary for the existence of stationary solutions of such random recurrence equations.

(b) The conditions of Theorem 3.1 imply the conditions (3.6) with the matrix norm defined as the natural norm $\|A\|_{B,r} = \|S^{-1}AS\|_r$, corresponding to the vector norm

$$\|\mathbf{c}\|_{B,r} := \|S^{-1}\mathbf{c}\|_r, \qquad \mathbf{c} \in \mathbb{C}^q. \tag{3.7}$$

Observe, however, that the conditions of Theorem 3.1 are in general not necessary for stationarity. For example, using methods similar to those in the proofs of Theorems 3.1 and 3.3, it can be shown that for any vector norm $\|\cdot\|$ and for $t \geq 0$,

$$\|J_{0,t}\| \leq \|e^{Bt}\| + e^{\|B\|t} \exp\left(\sum_{0<s\leq t} \log(1 + (\Delta L_s)^2 \|\mathbf{e}\mathbf{a}'\|)\right)$$
$$\times \|\mathbf{e}\mathbf{a}'\| \sum_{0<s\leq t} (\Delta L_s)^2.$$

Now if $\lambda(B) < 0$, then $\|e^{Bt}\| \to 0$ as $t \to \infty$, and (3.6) can be satisfied without assuming that all the eigenvalues of $B$ are distinct, but choosing $\|\mathbf{a}\|$ sufficiently small and imposing certain integrability conditions on $L$. We shall not pursue this argument here because the conditions of Theorem 3.1 will be sufficient for our purposes.

The matrices $J_{s,t}$ and the vector $\mathbf{K}_{s,t}$ of Theorem 3.3 will be constructed explicitly when $L$ is compound-Poisson, and in the general case will be obtained as the limit of the corresponding quantities for compound-Poisson-driven processes. In the compound-Poisson case we shall show that the stationary state vector satisfies a distributional fixed point equation which is much easier to handle than (3.5). Also, we compare the stationary distribution of $\mathbf{Y}_\infty$ with the stationary distribution of the state vector when sampled at the jump times of the Lévy process. This is the content of the next theorem.

THEOREM 3.5. *(a) Let $(\mathbf{Y}_t)_{t\geq 0}$ be the state process of a* COGARCH$(p,q)$ *process with parameters $B$, $\mathbf{a}$ and $\alpha_0$. Suppose that the Lévy measure $\nu_L$ of the driving Lévy process $L$ is finite and write the compound-Poisson process $[L,L]^{(d)}$ in the form*

$$[L,L]_t^{(d)} = \sum_{0<s\leq t} (\Delta L_s)^2 = \sum_{i=1}^{N(t)} Z_i,$$

where $N(t)$ is the number of jumps of $L$ in the time interval $(0,t]$ and $Z_i$ is the square of the $i$th jump size. Let $T_1$ denote the time at which the first jump occurs and let $T_j$, $j = 2, 3, \ldots$, be the time intervals between the $(j-1)$st and $j$th jumps. Furthermore, let $(T_0, Z_0)$ be independent of $(T_i, Z_i)_{i \in \mathbb{N}}$ with the same distribution as $(T_1, Z_1)$. For $i \in \mathbb{N}_0$, let

$$C_i = (I + Z_i \mathbf{e} \mathbf{a}') e^{BT_i},$$
$$\mathbf{D}_i = \alpha_0 Z_i \mathbf{e}$$

and $\Gamma_n = \sum_{i=1}^n T_i$ (where $\Gamma_0 := 0$). Then the discrete time process $(\mathbf{Y}_{\Gamma_n})_{n \in \mathbb{N}_0}$ satisfies the random recurrence equation

(3.8) $$\mathbf{Y}_{\Gamma_{n+1}} = C_{n+1} \mathbf{Y}_{\Gamma_n} + \mathbf{D}_{n+1}, \qquad n \in \mathbb{N}_0.$$

Furthermore, for any $t > 0$,

$$\mathbf{Y}_t = e^{B(t - \Gamma_{N(t)})} \left[ \mathbf{1}_{\{N(t) \neq 0\}} \mathbf{D}_{N(t)} + \sum_{i=0}^{N(t)-2} C_{N(t)} \cdots C_{N(t)-i} \mathbf{D}_{N(t)-i-1} \right.$$

(3.9)
$$\left. + C_{N(t)} \cdots C_1 \mathbf{Y}_0 \right]$$

$$\stackrel{d}{=} e^{B(t - \Gamma_{N(t)})} \left[ \mathbf{1}_{\{N(t) \neq 0\}} \mathbf{D}_1 + \sum_{i=1}^{N(t)-1} C_1 \cdots C_i \mathbf{D}_{i+1} + C_1 \cdots C_{N(t)} \mathbf{Y}_0 \right].$$

(b) *Assume additionally that the conditions of Theorem* 3.1 *are satisfied. Then the infinite sum* $\sum_{i=0}^\infty C_1 \cdots C_i \mathbf{D}_{i+1}$ *converges almost surely absolutely to a random vector* $\widehat{\mathbf{Y}}$, *which has the stationary distribution of the sequence* $(\mathbf{Y}_{\Gamma_n})_{n \in \mathbb{N}_0}$. *The stationary state vector* $\mathbf{Y}_\infty$ *satisfies*

(3.10) $$\mathbf{Y}_\infty \stackrel{d}{=} e^{BT} \widehat{\mathbf{Y}},$$

*where* $T$ *is independent of* $(T_i, Z_i)_{i \in \mathbb{N}_0}$ *and has the distribution of* $T_1$. *Furthermore,* $\mathbf{Y}_\infty$ *is the unique solution in distribution of the distributional fixed point equation*

(3.11) $$\mathbf{Y}_\infty \stackrel{d}{=} Q \mathbf{Y}_\infty + \mathbf{R},$$

*where* $\mathbf{Y}_\infty$ *is independent of* $(Q, \mathbf{R})$ *and*

$$Q := e^{BT_0}(I + Z_0 \mathbf{e} \mathbf{a}'),$$
$$\mathbf{R} := \alpha_0 Z_0 e^{BT_0} \mathbf{e}.$$

The fixed point equation (3.11) will play a crucial role in the determination of the covariance matrix of $\mathbf{Y}_\infty$, which is studied in the next section.

**4. Second-order properties of the volatility process.** In this section $(\mathbf{Y}_t)_{t\geq 0}$ denotes the state process defined by (2.3), with parameters $B$, $\mathbf{a}$ and $\alpha_0$, and driving Lévy process $L$ with Lévy measure $\nu_L$. The aim of this section is to study the autocorrelation function of the volatility process $(V_t)_{t\geq 0}$. We shall write

$$\mu := \int_{\mathbb{R}} y^2 \, d\nu_L(y) \quad \text{and} \quad \rho := \int_{\mathbb{R}} y^4 \, d\nu_L(y),$$

and, if $\mu < \infty$ (i.e., $EL_1^2 < \infty$),

(4.1) $$\widetilde{B} := B + \mu \mathbf{e}\mathbf{a}'.$$

Observe that $\widetilde{B}$ has the same form as $B$, but with last row given by $(-\beta_q + \mu\alpha_1, \ldots, -\beta_1 + \mu\alpha_q)$. We first give sufficient conditions for the moments of $\mathbf{Y}_t$ to exist.

PROPOSITION 4.1. *Suppose that the eigenvalues of $B$ are distinct, $\lambda = \lambda(B) < 0$, $\|\cdot\|$ is any vector norm on $\mathbb{C}^q$ and $k \in \mathbb{N}$. Then the following results hold.*

(a) *If $E|L_1|^{2k} < \infty$ and $E\|\mathbf{Y}_0\|^k < \infty$,*

$$E\|\mathbf{Y}_t\|^k < \infty \qquad \forall t \geq 0.$$

(b) *If $E|L_1|^{2k} < \infty$, $r \in [1,\infty]$, $S$ is a matrix such that $S^{-1}BS$ is diagonal and*

$$\int_{\mathbb{R}} ((1 + \|S^{-1}\mathbf{e}\mathbf{a}'S\|_r y^2)^k - 1) \, d\nu_L(y) < -\lambda k,$$

*then $S$ and $r$ satisfy (3.2) and $E\|\mathbf{Y}_\infty\|^k < \infty$. In particular, $E(\mathbf{Y}_\infty)$ exists if*

(4.2) $$EL_1^2 < \infty \quad \text{and} \quad \|S^{-1}\mathbf{e}\mathbf{a}'S\|_r \mu < -\lambda,$$

*and the covariance matrix $\mathrm{cov}(\mathbf{Y}_\infty)$ exists if*

(4.3) $$EL_1^4 < \infty \quad \text{and} \quad \|S^{-1}\mathbf{e}\mathbf{a}'S\|_r^2 \rho < 2(-\lambda - \|S^{-1}\mathbf{e}\mathbf{a}'S\|_r \mu).$$

*Furthermore, (4.3) implies (4.2), and (4.2) implies that all the eigenvalues of $\widetilde{B}$ have strictly negative real parts, in particular that $\widetilde{B}$ is invertible and $\beta_q \neq \alpha_1 \mu$.*

Next, we determine the autocovariance function of the (not necessarily stationary) volatility process of Definition 2.1.

THEOREM 4.2.  *Let $(V_t)_{t\geq 0}$ be the volatility process specified in Definition 2.1, with state process $(\mathbf{Y}_t)_{t\geq 0}$ and parameters $B$, $\mathbf{a}$ and $\alpha_0$. Suppose that $EL_1^4 < \infty$ and that $E\|\mathbf{Y}_t\|^2 < \infty \ \forall t \geq 0$ (as is the case, e.g., if the conditions of Proposition 4.1 are satisfied). Then, with $\widetilde{B}$ defined as in (4.1),*

$$\text{(4.4)} \qquad \operatorname{cov}(V_{t+h}, V_t) = \mathbf{a}' e^{\widetilde{B}h} \operatorname{cov}(\mathbf{Y}_t)\mathbf{a}, \qquad t, h \geq 0.$$

Since we are primarily interested in the stationary volatility process, we need to evaluate $\operatorname{cov}(\mathbf{Y}_\infty)$. However, first we need an expression for $E(\mathbf{Y}_\infty)$.

LEMMA 4.3.  *Suppose that all the eigenvalues of $B$ are distinct and that (4.2) holds. Then*

$$\text{(4.5)} \qquad E(\mathbf{Y}_\infty) = -\alpha_0 \mu \widetilde{B}^{-1} \mathbf{e} = \frac{\alpha_0 \mu}{\beta_q - \alpha_1 \mu} \mathbf{e}_1.$$

The following theorem contains the main results of this section. It demonstrates that the autocorrelation function of the stationary COGARCH volatility process is the same as that of a continuous-time ARMA process. This reflects the corresponding discrete-time result that the autocorrelation function of a GARCH volatility process is the same as that of a discrete-time ARMA process.

THEOREM 4.4.  *Suppose that the eigenvalues of the matrix $B$ are distinct, $\lambda(B) < 0$ and (4.3) holds. Then the matrix $(I \otimes \widetilde{B}) + (\widetilde{B} \otimes I) + \rho((\mathbf{ea}') \otimes (\mathbf{ea}'))$ is invertible and the covariance matrix of $\mathbf{Y}_\infty$ is the unique solution of*

$$\text{(4.6)} \quad \begin{aligned} &[(I \otimes \widetilde{B}) + (\widetilde{B} \otimes I) + \rho((\mathbf{ea}') \otimes (\mathbf{ea}'))]\operatorname{vec}(\operatorname{cov}(\mathbf{Y}_\infty)) \\ &= \frac{-\alpha_0^2 \beta_q^2 \rho}{(\beta_q - \mu\alpha_1)^2} \operatorname{vec}(\mathbf{ee}'). \end{aligned}$$

*Let $(\psi_t)_{t\geq 0}$ be a stationary $\operatorname{CARMA}(q, p-1)$ process (as defined in Section 2) with location parameter $0$, moving average coefficients $\alpha_1, \ldots, \alpha_p$, autoregressive coefficients $\beta_1 - \mu\alpha_q, \beta_2 - \mu\alpha_{q-1}, \ldots, \beta_q - \alpha_1\mu$, driving Lévy process $\widetilde{L}$ and corresponding state process $(\zeta_t)_{t\geq 0}$. Suppose that $E(\widetilde{L}_1)^2 < \infty$, $E(\widetilde{L}_1) = \mu$ and $\operatorname{var}(\widetilde{L}_1) = \rho$, and define*

$$m := \rho \int_0^\infty \mathbf{a}' e^{\widetilde{B}t} \mathbf{ee}' e^{\widetilde{B}'t} \mathbf{a}\, dt = \operatorname{var}(\psi_t).$$

*Then $0 \leq m < 1$ and*

$$\operatorname{cov}(\mathbf{Y}_\infty) = \frac{\alpha_0^2 \beta_q^2}{(\beta_q - \mu\alpha_1)^2(1-m)} \operatorname{cov}(\zeta_\infty)$$

$$\text{(4.7)} \quad = \frac{\alpha_0^2 \beta_q^2 \rho}{(\beta_q - \mu\alpha_1)^2(1-m)} \int_0^\infty e^{\widetilde{B}t} \mathbf{e}\mathbf{e}' e^{\widetilde{B}'t}\,dt,$$

$$\text{(4.8)} \quad \operatorname{var}(V_\infty) = \frac{\alpha_0^2 \beta_q^2}{(\beta_q - \mu\alpha_1)^2}\frac{m}{1-m},$$

$$\text{(4.9)} \quad E(V_\infty) = \frac{\alpha_0 \beta_q}{\beta_q - \mu\alpha_1},$$

$$\text{(4.10)} \quad E(\psi_\infty) = \frac{\alpha_1 \mu}{\beta_q - \mu\alpha_1}.$$

If $(V_t)_{t\geq 0}$ is the stationary COGARCH volatility process, then

$$\text{(4.11)} \quad \operatorname{cov}(V_{t+h}, V_t) = \frac{\alpha_0^2 \beta_q^2}{(\beta_q - \mu\alpha_1)^2(1-m)}\operatorname{cov}(\psi_{t+h}, \psi_t), \qquad t, h \geq 0,$$

showing, in particular, that $V$ has the same autocorrelation function as $\psi$. If the eigenvalues $\widetilde{\lambda}_1, \ldots, \widetilde{\lambda}_q$ of $\widetilde{B}$ are also distinct, and $a(z)$ and $\widetilde{b}(z)$ are the characteristic polynomials associated with $\mathbf{a}$ and $\widetilde{B}$, then

$$\text{(4.12)} \quad \begin{aligned} \operatorname{cov}(V_{t+h}, V_t) &= \frac{\alpha_0^2 \beta_q^2 \rho}{(\beta_q - \mu\alpha_1)^2(1-m)} \\ &\quad \times \sum_{j=1}^q \frac{a(\widetilde{\lambda}_j)a(-\widetilde{\lambda}_j)}{\widetilde{b}'(\widetilde{\lambda}_j)\widetilde{b}(-\widetilde{\lambda}_j)} e^{\widetilde{\lambda}_j h}, \qquad t, h \geq 0, \end{aligned}$$

where $\widetilde{b}'$ denotes the derivative of $\widetilde{b}$.

**5. Positivity conditions for the volatility.** For the definition of the COGARCH price process $dG_t = \sqrt{V_t}\,dt$ to make sense, it is necessary that $V_t$ be nonnegative for all $t \geq 0$. The following theorem gives necessary and sufficient conditions for this to occur with probability 1.

THEOREM 5.1. (a) *Let $(\mathbf{Y}_t)_{t\geq 0}$ be the state vector of a $\operatorname{COGARCH}(p,q)$ volatility process $(V_t)_{t\geq 0}$ with parameters $B$, $\mathbf{a}$ and $\alpha_0 > 0$. Let $\gamma \geq -\alpha_0$ be a real constant. Suppose that the following two conditions hold:*

$$\text{(5.1)} \quad \mathbf{a}' e^{Bt}\mathbf{e} \geq 0 \qquad \forall t \geq 0,$$

$$\text{(5.2)} \quad \mathbf{a}' e^{Bt}\mathbf{Y}_0 \geq \gamma \qquad \text{a.s. } \forall t \geq 0.$$

*Then for any driving Lévy process, with probability 1,*

$$\text{(5.3)} \quad V_t \geq \alpha_0 + \gamma \geq 0 \qquad \forall t \geq 0.$$

*Conversely, if either* (5.2) *fails or* (5.2) *holds with* $\gamma > -\alpha_0$ *and* (5.1) *fails, then there exists a driving compound-Poisson process* $L$ *and* $t_0 \geq 0$ *such that* $P(V_{t_0} < 0) > 0$.

(b) *Suppose that all the eigenvalues of* $B$ *are distinct and that* (3.2) *and* (5.1) *both hold. Then with probability* 1 *the stationary* COGARCH$(p,q)$ *volatility process* $(V_t)_{t \geq 0}$ *satisfies*

$$V_t \geq \alpha_0 > 0 \qquad \forall t \geq 0.$$

For the stationary COGARCH volatility process or for the process with $\mathbf{Y}_0 = \mathbf{0}$, the condition (5.1) alone is sufficient for almost sure nonnegativity. The expression $\mathbf{a}' e^{Bt} \mathbf{e}$ is in fact the kernel of a CARMA process with autoregressive coefficients $b_1, \ldots, b_q$ and moving average coefficients $a_1, \ldots, a_q$. Results that pertain to nonnegativity of a CARMA kernel were recently obtained by Tsai and Chan [28]. We state their results in the next theorem in the context of COGARCH rather than CARMA processes. Statement (e) below was also obtained by Todorov and Tauchen [27]. Recall that a function $\phi$ on $(0, \infty)$ is called *completely monotone* if it possesses derivatives of all orders and satisfies $(-1)^n (d^n \phi / dt^n)(t) \geq 0$ for all $t > 0$ and all $n \in \mathbb{N}_0$.

THEOREM 5.2. *Let* $B$ *and* $\mathbf{a}$ *be the parameters of a* COGARCH$(p,q)$ *process. If* $\lambda(B) < 0$, *and* $\alpha_1 > 0$, *we have the following results.*

(a) *For the* COGARCH$(p,q)$ *process,* (5.1) *holds if and only if the ratio of the characteristic polynomials* $a(\cdot)/b(\cdot)$ *is completely monotone on* $(0, \infty)$.

(b) *A sufficient condition for* (5.1) *to hold for the* COGARCH$(1,q)$ *process is that either* (i) *all eigenvalues of* $B$ *are real and negative or* (ii) *if* $(\lambda_{i_1}, \lambda_{i_1+1}), \ldots, (\lambda_{i_r}, \lambda_{i_r+1})$ *is a partition of the set of all pairs of complex conjugate eigenvalues of* $B$ *(counted with multiplicity), then there exists an injective mapping* $u : \{1, \ldots, r\} \to \{1, \ldots, q\}$ *such that* $\lambda_{u(j)}$ *is a real eigenvalue of* $B$ *satisfying* $\lambda_{u(j)} \geq \Re(\lambda_{i_j})$.

(c) *A necessary condition for* (5.1) *to hold for the* COGARCH$(1,q)$ *process is that there exists a real eigenvalue of* $B$ *not smaller than the real part of all other eigenvalues of* $B$.

(d) *Suppose* $2 \leq p \leq q$, *that all eigenvalues of* $B$ *are negative and ordered as in Definition* 2.3, *and that the roots* $\gamma_j$ *of* $a(z) = 0$ *are negative and ordered such that* $\gamma_{p-1} \leq \cdots \leq \gamma_1 < 0$. *Then a sufficient condition for* (5.1) *to hold for the* COGARCH$(p,q)$ *process is that*

$$\sum_{i=1}^{k} \gamma_i \leq \sum_{i=1}^{k} \lambda_i \qquad \forall k \in \{1, \ldots, p-1\}.$$

(e) *A necessary and sufficient condition for* (5.1) *in the* COGARCH$(2,2)$ *case is that both eigenvalues of* $B$ *are real, that* $\alpha_2 \geq 0$ *and that* $\alpha_1 \geq -\alpha_2 \lambda(B)$.

Although characterization (a) may be difficult to check in general, it gives a method for constructing further pairs $(\mathbf{a}, B)$ for which (5.1) holds, since the product of two completely monotone functions is again completely monotone.

**6. The autocorrelation of the squared increments.** In Section 4 we investigated the behavior of the autocorrelation function of the volatility process. Since one of the striking features of observed financial time series is that the returns have negligible correlation while the squared returns are significantly correlated, we now turn to the second-order properties of the increments of the COGARCH process itself. We therefore assume that $V$ is strictly stationary and nonnegative, and define, for $r > 0$,

$$G_t^{(r)} := G_{t+r} - G_t = \int_{(t,t+r]} \sqrt{V_s}\, dL_s, \qquad t \geq 0.$$

It is easy to see that $(G_t^{(r)})_{t \geq 0}$ is a stationary process. Let $\mu$ and $\widetilde{B}$ be defined as in Section 4. We then have the following theorem.

THEOREM 6.1. *Let $B$, $\mathbf{a}$ and $\alpha_0$ be the parameters of a* COGARCH$(p,q)$ *process whose driving Lévy process has mean zero. Suppose that the eigenvalues of $B$ are distinct, that (4.2) and (5.1) hold, and that $V$ is the stationary volatility process. Then for any $t \geq 0$ and $h \geq r > 0$,*

(6.1) $$E(G_t^{(r)}) = 0,$$

(6.2) $$E((G_t^{(r)})^2) = \frac{\alpha_0 \beta_q r}{\beta_q - \mu \alpha_1} E(L_1^2),$$

(6.3) $$\operatorname{cov}(G_t^{(r)}, G_{t+h}^{(r)}) = 0.$$

*If in addition* (4.3) *holds, then*

(6.4) $$\operatorname{cov}((G_t^{(r)})^2, (G_{t+h}^{(r)})^2) = \mathbf{a}' e^{\widetilde{B}h} \mathbf{H}_r, \qquad h \geq r,$$

*where*

$$\mathbf{H}_r := E(L_1^2)\, \widetilde{B}^{-1}(I - e^{-\widetilde{B}r}) \operatorname{cov}(\mathbf{Y}_r, G_r^2).$$

The autocovariance function (6.4), like that of the CARMA process with parameters $\widetilde{B}$ and $\mathbf{a}$, is a linear combination of terms of the form $e^{\tilde{\lambda}_j h}$, $j = 1, \ldots, q$, where $\tilde{\lambda}_1, \ldots, \tilde{\lambda}_q$ are the eigenvalues of $\widetilde{B}$.

**7. An example.** In this section we illustrate the properties established above using the COGARCH(1,3) process driven by a compound-Poisson process with jump rate 2 and normally distributed jumps with mean 0 and variance 0.74. The COGARCH coefficients are $\alpha_0 = \alpha_1 = 1, \beta_1 = 1.2, \beta_2 = 0.48 + \pi^2$ and $\beta_3 = 0.064 + 0.4\pi^2$, from which we find that the eigenvalues of $B$ are $-0.4$, $-0.4 + \pi i$ and $-0.4 - \pi i$. With $S$ defined as in (3.1), $\|S^{-1}\mathbf{ea}'S\|_2 = 0.21493$ and it is easy to check from this that the conditions (4.2) and (4.3) are satisfied. Condition (b)(ii) of Theorem 5.2 also implies that the volatility process is nonnegative.

The eigenvalues of the matrix $\widetilde{B} = B + \mu \mathbf{ea}'$ are $-0.25038$, $-0.47481 + 3.14426i$ and $-0.47481 - 3.14426i$. From (4.12) we conclude that the autocorrelation of the volatility in this case is a linear combination of $\exp(-0.25038t)$, and a damped sinusoid with period approximately equal to 2 and damping factor $\exp(-0.47481t)$.

The top graph in Figure 1 shows the values at integer times $101, \ldots, 8100$ of a simulated series $(G_t)$ with the parameters specified above, $\mathbf{Y}_0 = (1, 1, 1)'$ and $G(0) = 0$. The second graph shows the differenced series $(G_{t+1} - G_t)_{t=100,\ldots,8099}$ and the last graph shows the volatility $(\sigma_t^2)_{t=101,\ldots,8100}$.

As is the case for a discrete-time GARCH process, the increments $(G_{t+1} - G_t)$ exhibit no significant correlation, but the squared increments $((G_{t+1} - G_t)^2)$ have highly significant correlations as shown in the second graph of Figure 2. The first graph in Figure 2 shows the sample autocorrelation function of the volatility process at integer lags. This too is highly significant for large lags, reflecting the long-memory property frequently observed in financial time series. As expected from the remarks in the first paragraph above, it has the form of an exponentially decaying term plus a small damped sinusoidal term with period approximately equal to 2.

**8. Proofs for Section 3.** We start by proving Theorem 3.5, since (3.9) will be needed in the proof of Theorems 3.1 and 3.3.

PROOF OF THEOREM 3.5. (a) It follows from (2.3) that $\mathbf{Y}_t$ satisfies $d\mathbf{Y}_t = B\mathbf{Y}_t \, dt$ for $t \in [\Gamma_n, \Gamma_{n+1})$, so that

$$(8.1) \qquad \mathbf{Y}_t = e^{B(t-\Gamma_n)} \mathbf{Y}_{\Gamma_n}, \qquad t \in [\Gamma_n, \Gamma_{n+1}), n \in \mathbb{N}_0.$$

At time $\Gamma_{n+1}$ a jump of size $\mathbf{e}(\alpha_0 + \mathbf{a}'\mathbf{Y}_{\Gamma_{n+1}-})Z_{n+1}$ occurs, so that

$$\begin{aligned}
\mathbf{Y}_{\Gamma_{n+1}} &= \mathbf{Y}_{\Gamma_{n+1}-} + \mathbf{e}(\alpha_0 + \mathbf{a}'\mathbf{Y}_{\Gamma_{n+1}-})Z_{n+1} \\
&= (I + Z_{n+1}\mathbf{ea}')\mathbf{Y}_{\Gamma_{n+1}-} + \alpha_0 Z_{n+1}\mathbf{e} \\
&= C_{n+1}\mathbf{Y}_{\Gamma_n} + \mathbf{D}_{n+1}, \qquad n \in \mathbb{N}_0,
\end{aligned}$$

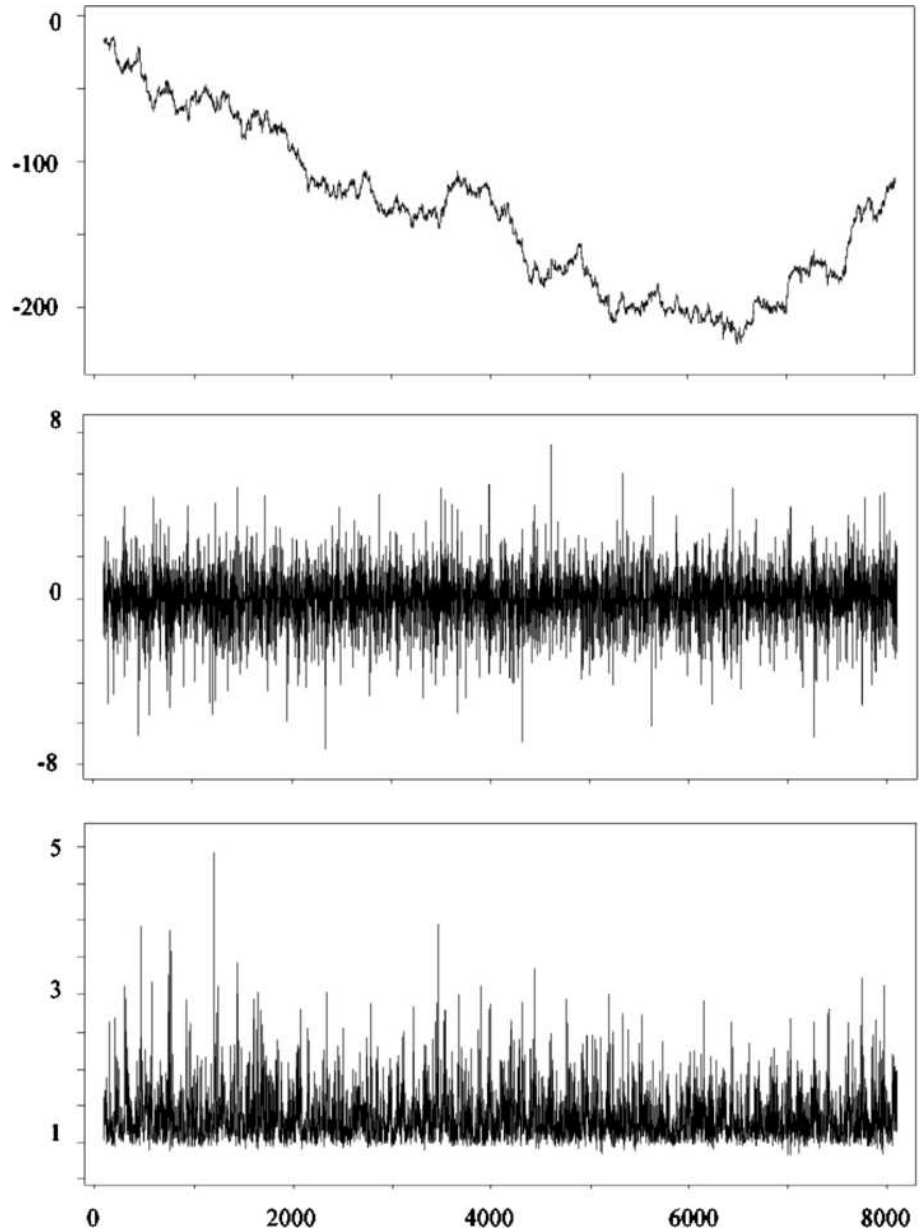

Fig. 1. *The simulated compound-Poisson driven* COGARCH(1,3) *process with jump rate* 2, *normally distributed jumps with mean* 0 *and variance* 0.74, *and coefficients* $\alpha_0 = \alpha_1 = 1, \beta_1 = 1.2, \beta_2 = 0.48 + \pi^2$ *and* $\beta_3 = 0.064 + 0.4\pi^2$. *The graphs show the process* $(G_t)$ *sampled at integer times* (top), *the corresponding increments* $((G_{t+1} - G_t))$ (center) *and the corresponding volatility sequence* $(V_t = \sigma_t^2)$ (bottom).

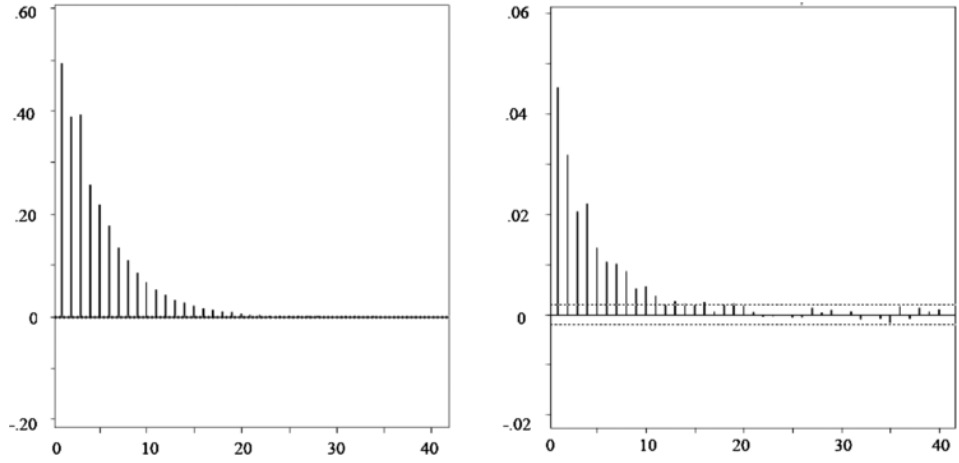

Fig. 2. *The sample autocorrelation functions of the volatilities $(V_t)$ (left) and of the squared COGARCH increments $((G_{t+1} - G_t)^2)$ (right) of a realization of length* 1,000,000 *of the COGARCH process with parameters as specified in Figure* 1.

which is (3.8). Solving this recursion gives

$$\mathbf{Y}_{\Gamma_n} = \mathbf{D}_n + \sum_{i=0}^{n-2} C_n \cdots C_{n-i} \mathbf{D}_{n-i-1} + C_n \cdots C_1 \mathbf{Y}_0, \qquad n \in \mathbb{N},$$

and the first equality in (3.9) follows from this and $\mathbf{Y}_t = e^{B(t - \Gamma_{N(t)})} \mathbf{Y}_{\Gamma_{N(t)}}$. The second equality in (3.9) is a consequence of the fact that the infinite random element $(N(t), \Gamma_{N(t)}, C_{N(t)}, \mathbf{D}_{N(t)}, \ldots, C_1, \mathbf{D}_1, 0, 0, \ldots)$ has the same distribution as $(N(t), \Gamma_{N(t)}, C_1, \mathbf{D}_1, \ldots, C_{N(t)}, \mathbf{D}_{N(t)}, 0, 0, \ldots)$; indeed, for any $n \in \mathbb{N}_0$ and $c \geq 0$, the random vectors $(C_1, \mathbf{D}_1), \ldots, (C_n, \mathbf{D}_n)$ are i.i.d. and depend on the restriction $\{N(t) = n, \ \Gamma_{N(t)} \geq c\}$ only in terms of $\sum_{i=1}^n T_i$ and $T_{n+1}$, but not on the $T_i$, $i = 1, \ldots, n$, individually.

(b) Let $S$ be such that $S^{-1}BS =: \Lambda$ is diagonal and define the vector norm $\|\mathbf{c}\|_{B,r} = \|S^{-1}\mathbf{c}\|_r$ as in (3.7), so that the associated natural matrix norm is $\|A\|_{B,r} = \|S^{-1}AS\|_r$. Then we have for $t \geq 0$,

(8.2) $$\|e^{Bt}\|_{B,r} = \|Se^{\Lambda t}S^{-1}\|_{B,r} = \|e^{\Lambda t}\|_r = e^{\lambda t}.$$

This gives $\|C_1\|_{B,r} \leq (1 + Z_1\|\mathbf{ea}'\|_{B,r})e^{\lambda T_1}$ and $\|\mathbf{D}_1\|_{B,r} = \alpha_0\|\mathbf{e}\|_{B,r} Z_1$, so that, using $\nu_{[L,L]}([x, \infty)) = \nu_L\{y \in \mathbb{R} : |y| \geq \sqrt{x}\}$ for $x \geq 0$,

$$E \log \|C_1\|_{B,r} \leq \lambda E(T_1) + E \log(1 + Z_1\|\mathbf{ea}'\|_{B,r})$$
$$= \frac{\lambda}{\nu_L(\mathbb{R})} + \frac{1}{\nu_L(\mathbb{R})} \int_{(0,\infty)} \log(1 + \|\mathbf{ea}'\|_{B,r}\, y^2)\, d\nu_L(y) < 0$$

by (3.2) and

$$E\log^+(Z_1) = \frac{1}{\nu_L(\mathbb{R})} \int_{\mathbb{R}} \log^+(y^2)\,d\nu_L(y) < \infty.$$

From the general theory of random recurrence equations, this implies the almost sure absolute convergence of $\sum_{i=0}^{\infty} C_1 \cdots C_i \mathbf{D}_{i+1}$ to $\widehat{\mathbf{Y}}$, which has the stationary distribution of $(\mathbf{Y}_{\Gamma_n})_{n\in\mathbb{N}}$ (see, e.g., [7]).

To prove (3.10), for $m \in \mathbb{N}$, let

$$\widehat{\mathbf{Y}}_m := \sum_{i=0}^{m-1} C_1 \cdots C_i D_{i+1} + C_1 \cdots C_m \mathbf{Y}_0$$

and

$$\mathbf{Y}_{t,m} := e^{B(t-\Gamma_{N(t)})}\widehat{\mathbf{Y}}_m, \qquad t \geq 0.$$

Since the random variable $(t - \Gamma_{N(t)})$ is asymptotically independent of $T_1$, $Z_1, \ldots, T_m, Z_m$ (for $t \to \infty$, $m$ fixed), it follows that $e^{B(t-\Gamma_{N(t)})}$ is asymptotically independent of $\widehat{\mathbf{Y}}_m$ and, hence, $\mathbf{Y}_{t,m}$ converges in distribution to $e^{BT}\widehat{\mathbf{Y}}_m$ as $t \to \infty$, where $T$ is exponentially distributed with parameter $\nu_L(\mathbb{R})$ (e.g., [26], Section 7.4.4) and is independent of $T_1, Z_1, \ldots, T_m, Z_m$ and, hence, can be chosen to be independent of $(T_i)_{i\in\mathbb{N}}$, $(Y_i)_{i\in\mathbb{N}}$ (as in the statement of the theorem). Moreover, $e^{BT}\widehat{\mathbf{Y}}_m$ converges almost surely, hence in distribution to $e^{BT}\widehat{\mathbf{Y}}$, as $m \to \infty$. Denote by $\widetilde{\mathbf{Y}}_t$ the expression in the lower line of (3.9). Then (3.10) and, in particular, the existence of the limit variable $\mathbf{Y}_{\infty}$ in the compound-Poisson case follow from (3.9) and a variant of Slutsky's theorem (e.g., [11], Proposition 6.3.9), provided

(8.3) $$\lim_{m\to\infty} \limsup_{t\to\infty} P(\|\widetilde{\mathbf{Y}}_t - \mathbf{Y}_{t,m}\|_{B,r} > \varepsilon) = 0 \qquad \forall \varepsilon > 0.$$

Since $\|e^{B(t-\Gamma_{N(t)})}\|_{B,r} \leq 1$ and $\mathbf{1}_{\{N(t)\neq 0\}}\mathbf{D}_1 + \sum_{i=1}^{N(t)-1} C_1 \cdots C_i \mathbf{D}_{i+1} + C_1 \cdots C_{N(t)}\mathbf{Y}_0 - \widehat{\mathbf{Y}}_m$ converges almost surely, hence in probability as $t \to \infty$ to $\sum_{i=m}^{\infty} C_1 \cdots C_i \mathbf{D}_{i+1} - C_1 \cdots C_m \mathbf{Y}_0$, which itself converges almost surely to 0 as $m \to \infty$, (8.3) is true and (3.10) follows. That $\mathbf{Y}_{\infty}$ satisfies (3.11) is clear from (3.10); that it is the unique solution follows from $E\log\|Q\|_{B,r} < 0$ and $E\log^+\|\mathbf{R}\|_{B,r} < \infty$. □

The proof of Theorem 3.5(b) already showed the existence of the limit variable $\mathbf{Y}_{\infty}$ for the case of a driving compound-Poisson process. Nevertheless, this existence will be reestablished in the proof of Theorems 3.1 and 3.3 for the general case, making use of Theorem 3.5(a) only. We shall use an approximation argument and introduce the following notation:

DEFINITION 8.1. Let $L$ be a Lévy process. Then for any $\varepsilon > 0$, the $\sqrt{\varepsilon}$-cut Lévy process $(L_t^{(\varepsilon)})_{t \geq 0}$ is defined by

$$L_t^{(\varepsilon)} := \sum_{0 < s \leq t, |\Delta L_s| \geq \sqrt{\varepsilon}} \Delta L_s, \qquad t \geq 0.$$

If $(\mathbf{Y}_t)_{t \geq 0}$ is a state process of a COGARCH$(p, q)$ process driven by $L$, then the state process of the COGARCH$(p, q)$ process with the same parameters and starting vector but driving Lévy process $(L_t^{(\varepsilon)})_{t \geq 0}$ will be denoted by $(\mathbf{Y}_t^{(\varepsilon)})_{t \geq 0}$.

The quadratic covariation of $L^{(\varepsilon)}$ is given by

$$[L^{(\varepsilon)}, L^{(\varepsilon)}]_t = [L^{(\varepsilon)}, L^{(\varepsilon)}]_t^{(d)}$$
$$= \sum_{0 < s \leq t, |\Delta L_s|^2 \geq \varepsilon} |\Delta L_s|^2.$$

In particular, the corresponding COGARCH volatility is driven by a compound-Poisson process. With this notation, we have the following lemma:

LEMMA 8.2. *Let $(\mathbf{Y}_t)_{t \geq 0}$ be the state process of a COGARCH$(p, q)$ process. Then $\mathbf{Y}_t^{(\varepsilon)}$ converges in ucp to $\mathbf{Y}_t$ as $\varepsilon \to 0$.*

PROOF. This is an easy consequence of perturbation results in stochastic differential equations: recalling the definition of prelocal convergence in $\underline{H}^p$, $1 \leq p < \infty$, as in [24], page 260, it is easy to see that $[L^{(\varepsilon)}, L^{(\varepsilon)}]$ converges prelocally to $[L, L]^{(d)}$ in $\underline{H}^p$, $1 \leq p < \infty$, as $\varepsilon \to 0$ (e.g., with stopping times $T_k = k$). The claim then follows from Theorems 14 and 15 of Chapter V in [24]. □

PROOF OF THEOREMS 3.1 AND 3.3. We shall first concentrate on (3.3) and (3.4), and then prove Theorem 3.1 and the rest of Theorem 3.3 simultaneously. Let $\varepsilon > 0$ and assume the representation

$$[L^{(\varepsilon)}, L^{(\varepsilon)}]_t = \sum_{i=1}^{N_\varepsilon(t)} Z_i^{(\varepsilon)},$$

where $L^{(\varepsilon)}$ is the $\sqrt{\varepsilon}$-cut Lévy process of Definition 8.1. Define $C_i^{(\varepsilon)}$ and $\mathbf{D}_i^{(\varepsilon)}$ similarly as in Theorem 3.5. Furthermore, let

$$J_{0,t}^{(\varepsilon)} := e^{B(t - \Gamma_{N_\varepsilon(t)}^{(\varepsilon)})} C_{N_\varepsilon(t)}^{(\varepsilon)} \cdots C_1^{(\varepsilon)},$$

$$\mathbf{K}_{0,t}^{(\varepsilon)} := e^{B(t - \Gamma_{N_\varepsilon(t)}^{(\varepsilon)})} \bigg[ \mathbf{1}_{\{N_\varepsilon(t) \neq 0\}} \mathbf{D}_{N_\varepsilon(t)}^{(\varepsilon)}$$

$$+ \sum_{i=0}^{N_\varepsilon(t)-2} C^{(\varepsilon)}_{N_\varepsilon(t)} \cdots C^{(\varepsilon)}_{N_\varepsilon(t)-i} \mathbf{D}^{(\varepsilon)}_{N_\varepsilon(t)-i-1} \Bigg].$$

Then, by Theorem 3.5(a),

(8.4) $$\mathbf{Y}^{(\varepsilon)}_t = J^{(\varepsilon)}_{0,t} \mathbf{Y}_0 + \mathbf{K}^{(\varepsilon)}_{0,t}.$$

From the previous lemma we know that $\mathbf{Y}^{(\varepsilon)}_t$ converges in ucp to $\mathbf{Y}_t$ as $\varepsilon \to 0$. Since this is true for any starting value $\mathbf{Y}_0$, it holds in particular for $\mathbf{Y}_0 = 0$, and from (8.4) it follows that $\mathbf{K}^{(\varepsilon)}_{0,t}$ converges in ucp to some $\mathbf{K}_{0,t}$ as $\varepsilon \to 0$. Hence, again from (8.4), it follows that for arbitrary $\mathbf{Y}_0$,

$$J^{(\varepsilon)}_{0,t} \mathbf{Y}_0 = \mathbf{Y}^{(\varepsilon)}_t - \mathbf{K}^{(\varepsilon)}_{0,t} \overset{ucp}{\to} \mathbf{Y}_t - \mathbf{K}_{0,t} \qquad \text{as } \varepsilon \to 0.$$

Since this holds for arbitrary $\mathbf{Y}_0$, we conclude that $J^{(\varepsilon)}_{0,t}$ converges in ucp to some $J_{0,t}$ as $\varepsilon \to 0$. From (8.4) it then follows that

$$\mathbf{Y}_t = J_{0,t} \mathbf{Y}_0 + \mathbf{K}_{0,t}.$$

By starting at an arbitrary time $s$ instead of at time 0, we obtain (3.3). For example, $J^{(\varepsilon)}_{s,t}$ is given by

$$J^{(\varepsilon)}_{s,t} = e^{B(t-\Gamma^{(\varepsilon)}_{N_\varepsilon(t)})} C_{N_\varepsilon(t)} \cdots C_{N_\varepsilon(s)+2}(I + Z_{N_\varepsilon(s)+1} \mathbf{ea}') e^{B(\Gamma^{(\varepsilon)}_{N_\varepsilon(s)+1}-s)},$$
$$0 \le s \le t,$$

giving (3.4). The independence and stationarity assertions on $(J_{s,t}, \mathbf{K}_{s,t})$ are clear, since $J_{s,t}$ and $\mathbf{K}_{s,t}$ are constructed only from the segment $(L_u)_{s<u\le t}$ of the Lévy process $L$.

Now assume that all eigenvalues of $B$ are distinct and that (3.2) holds. Applying (8.2) to $J^{(\varepsilon)}_{0,t}$ gives

$$\|J^{(\varepsilon)}_{0,t}\|_{B,r} \le \|e^{B(t-\Gamma^{(\varepsilon)}_{N_\varepsilon(t)})}\|_{B,r} \|C^{(\varepsilon)}_{N_\varepsilon(t)}\|_{B,r} \cdots \|C^{(\varepsilon)}_1\|_{B,r}$$

$$\le e^{\lambda(t-\Gamma^{(\varepsilon)}_{N_\varepsilon(t)})} \prod_{i=1}^{N_\varepsilon(t)} ((1+Z^{(\varepsilon)}_i \|\mathbf{ea}'\|_{B,r}) e^{\lambda(\Gamma^{(\varepsilon)}_i - \Gamma^{(\varepsilon)}_{i-1})})$$

(8.5) $$= e^{\lambda t} \exp\left(\sum_{i=1}^{N_\varepsilon(t)} \log(1 + Z^{(\varepsilon)}_i \|S^{-1}\mathbf{ea}'S\|_r)\right)$$

(8.6) $$\le e^{\lambda t} \exp\left(\sum_{0<s\le t} \log(1 + (\Delta L_s)^2 \|S^{-1}\mathbf{ea}'S\|_r)\right).$$

Since $\|J_{0,t}\|_{B,r} \leq \limsup_{\varepsilon \to 0} \|J_{0,t}^{(\varepsilon)}\|_{B,r}$, we conclude that

$$(8.7) \qquad \log \|J_{0,t}\|_{B,r} \leq \lambda t + \sum_{0 < s \leq t} \log(1 + (\Delta L_s)^2 \|S^{-1}\mathbf{e}\mathbf{a}'S\|_r),$$

giving

$$E \log \|J_{0,t}\|_{B,r} \leq t\left(\lambda + \int_{\mathbb{R}} \log(1 + \|S\mathbf{e}\mathbf{a}'S^{-1}\|y^2) \, d\nu_L(y)\right) < 0$$

by (3.2) (see, e.g., [24], Chapter I, Theorems 36 and 38). This is the left-hand inequality of (3.6). To show that $E \log^+ \|\mathbf{K}_{0,t}\|_{B,r} < \infty$, observe that

$$\|\mathbf{K}_{0,t}^{(\varepsilon)}\|_{B,r} \leq e^{\lambda(t - \Gamma_{N_\varepsilon(t)}^{(\varepsilon)})} \mathbf{1}_{\{N_\varepsilon(t) \neq 0\}} \alpha_0 \|\mathbf{e}\|_{B,r} Z_{N_\varepsilon(t)}^{(\varepsilon)}$$

$$(8.8) \qquad + \alpha_0 \|\mathbf{e}\|_{B,r} \sum_{i=0}^{N_\varepsilon(t)-2} e^{\lambda(t - \Gamma_{N_\varepsilon(t)-i-1}^{(\varepsilon)})} (1 + Z_{N_\varepsilon(t)}^{(\varepsilon)} \|\mathbf{e}\mathbf{a}'\|_{B,r})$$

$$\times \cdots \times (1 + Z_{N_\varepsilon(t)-i}^{(\varepsilon)} \|\mathbf{e}\mathbf{a}'\|_{B,r}) Z_{N_\varepsilon(t)-i-1}^{(\varepsilon)}$$

$$\leq \alpha_0 \|\mathbf{e}\|_{B,r} \mathbf{1}_{\{N_\varepsilon(t) \neq 0\}} Z_{N_\varepsilon(t)}^{(\varepsilon)}$$

$$+ \alpha_0 \|\mathbf{e}\|_{B,r} \sum_{i=0}^{N_\varepsilon(t)-2} \exp\left[\sum_{0 < s \leq t} \log(1 + (\Delta L_s)^2 \|\mathbf{e}\mathbf{a}'\|_{B,r})\right] Z_{N_\varepsilon(t)-i-1}^{(\varepsilon)}$$

$$\leq \alpha_0 \|S^{-1}\mathbf{e}\|_r \exp\left[\sum_{0 < s \leq t} \log(1 + (\Delta L_s)^2 \|S^{-1}\mathbf{e}\mathbf{a}'S\|_r)\right]$$

$$(8.9) \qquad \times \sum_{0 < s \leq t} (\Delta L_s)^2.$$

From this it follows that

$$\log \|\mathbf{K}_{0,t}\|_{B,r} \leq \log(\alpha_0 \|S^{-1}\mathbf{e}\|_r)$$

$$+ \sum_{0 < s \leq t} \log(1 + (\Delta L_s)^2 \|S^{-1}\mathbf{e}\mathbf{a}'S\|_r) + \log[L,L]_t^{(d)}.$$

The expectation of the second summand is finite as shown above and $E(\log[L,L]_t^{(d)}) < \infty$ since $\int_{(1,\infty)} \log x \, d\nu_{[L,L]}(x) = \int_{\mathbb{R}\setminus[-1,1]} \log x^2 \, d\nu_L(x) < \infty$, showing the right-hand inequality of (3.6).

Let $(J_n, \mathbf{K}_n)_{n \in \mathbb{N}}$ be an i.i.d. sequence with distribution $(J_{0,1}, \mathbf{K}_{0,1})$, independent of $L$ and $\mathbf{Y}_0$. Let $\gamma \in [0,1)$ and $n \in \mathbb{N}$. Then it follows from (3.3) that

$$\mathbf{Y}_{n+\gamma} = \mathbf{K}_{n+\gamma-1, n+\gamma} + \sum_{i=0}^{n-2} J_{n+\gamma-1, n+\gamma} \cdots J_{n+\gamma-i-1, n+\gamma-i} \mathbf{K}_{n+\gamma-i-2, n+\gamma-i-1}$$

$$+ J_{n+\gamma-1,n+\gamma} \cdots J_{\gamma,\gamma+1} \mathbf{Y}_\gamma$$

$$\stackrel{d}{=} \mathbf{K}_1 + \sum_{i=1}^{n-1} J_1 \cdots J_i \mathbf{K}_{i+1} + J_1 \cdots J_n \mathbf{Y}_\gamma$$

$$=: \mathbf{G}_n + H_n \mathbf{Y}_\gamma, \qquad \text{say.}$$

Since $E \log \|J_1\|_{B,r} < 0$ and $E \log^+ \|\mathbf{K}_1\|_{B,r} < \infty$, it follows from the general theory of random recurrence equations (e.g., [7]) that $H_n$ converges almost surely to 0 as $n \to \infty$ and that $\mathbf{G}_n$ converges almost surely absolutely to some random vector $\mathbf{G}$ as $n \to \infty$. Since $\mathbf{Y}$ has cadlag paths, it follows that $\sup_{\gamma \in [0,1)} \|\mathbf{Y}_\gamma\|_{B,r}$ is almost surely finite. Hence

$$\lim_{n \to \infty} \sup_{\gamma \in [0,1)} \|H_n \mathbf{Y}_\gamma\|_{B,r} = 0 \qquad \text{a.s.,}$$

and it follows that $\mathbf{Y}_t$ converges in distribution to $\mathbf{Y}_\infty := \mathbf{G}$ as $t \to \infty$. That $\mathbf{Y}_\infty$ satisfies (3.5) and is the unique solution is clear by the theory of random recurrence equations. Equations (3.5) and (3.3) then imply that if $\mathbf{Y}_0 \stackrel{d}{=} \mathbf{Y}_\infty$, then $\mathbf{Y}_t \stackrel{d}{=} \mathbf{Y}_\infty$ for all $t > 0$, showing strict stationarity of $(\mathbf{Y}_t)_{t \geq 0}$ since it is a Markov process. □

**9. Proofs for Section 4.** To prove Proposition 4.1, we will show that the state process $(\mathbf{Y}_t)_{t \geq 0}$ can be majorized by the state process of a COGARCH(1, 1) process, for which we can apply the moment conditions of Klüppelberg, Lindner and Maller [19]. We further show that under the conditions of Theorem 3.1, the stationary distribution $\mathbf{Y}_\infty$ can be approximated by stationary distributions of compound Poisson-driven COGARCH processes and that there is a majorant for this approximation. This will allow us to restrict attention to compound Poisson-driven processes when calculating autocorrelations, the general case following from Lebesgue's dominated convergence theorem. This is the content of the next lemma.

LEMMA 9.1. *Let $(\mathbf{Y}_t)_{t \geq 0}$ be the state process of a COGARCH$(p, q)$ process with parameters $B$, $\mathbf{a}$ and $\alpha_0 > 0$ such that all eigenvalues of $B$ are distinct and that $\lambda = \lambda(B) < 0$. Let $r \in [1, \infty]$ and let $S$ be such that $S^{-1}BS$ is diagonal. Denote by $\|\cdot\|_{B,r}$ the vector norm defined in (3.7). Furthermore, denote by $(\overline{\mathbf{Y}}_t)_{t \geq 0}$ the state process of a COGARCH$(1, 1)$ process satisfying (2.3) with the parameters $(B, \mathbf{a}, \alpha_0)$ replaced by $(\lambda, \|\mathbf{ea}'\|_{B,r}, \alpha_0 \|\mathbf{e}\|_{B,r})$ and initial state vector $\overline{\mathbf{Y}}_0 := \|\mathbf{Y}_0\|_{B,r}$. Then*

$$(9.1) \qquad \|\mathbf{Y}_t\|_{B,r} \leq \overline{\mathbf{Y}}_t, \qquad t \geq 0.$$

*If (3.2) is satisfied for this $r$, then there exist versions of $\mathbf{Y}_\infty$ and $\overline{\mathbf{Y}}_\infty$ such that*

$$(9.2) \qquad \|\mathbf{Y}_\infty\|_{B,r} \leq \overline{\mathbf{Y}}_\infty.$$

*Furthermore, if* $(\mathbf{Y}_t^{(\varepsilon)})_{t\geq 0}$ *is the process defined in Definition* 8.1 *for* $\varepsilon > 0$, *then versions of* $\mathbf{Y}_\infty^{(\varepsilon)}$ *can be chosen such that* $\|\mathbf{Y}_\infty^{(\varepsilon)}\|_{B,r} \leq \overline{\mathbf{Y}}_\infty$ *for all* $\varepsilon > 0$ *and* $\mathbf{Y}_\infty^{(\varepsilon)} \xrightarrow{P} \mathbf{Y}_\infty$, *as* $\varepsilon \to 0$.

PROOF. We use the notation and setup of the proof of Theorems 3.1 and 3.3. Let $\varepsilon > 0$ and define a COGARCH(1,1) state process $\overline{\mathbf{Y}}^{(\varepsilon)}$ similarly as above (with respect to $\mathbf{Y}^{(\varepsilon)}$). Let $\overline{J}_{0,t}^{(\varepsilon)}$ and $\overline{\mathbf{K}}_{0,t}^{(\varepsilon)}$ be defined similarly as $J_{0,t}^{(\varepsilon)}$ and $\mathbf{K}_{0,t}^{(\varepsilon)}$ (with respect to $\overline{\mathbf{Y}}^{(\varepsilon)}$). Then it is easy to see that $\overline{J}_{0,t}^{(\varepsilon)}$ and $\overline{\mathbf{K}}_{0,t}^{(\varepsilon)}$ are the right-hand sides of (8.5) and (8.8), respectively. In particular, $\|J_{0,t}^{(\varepsilon)}\|_{B,r} \leq \overline{J}_{0,t}^{(\varepsilon)}$ and $\|\mathbf{K}_{0,t}^{(\varepsilon)}\|_{B,r} \leq \overline{\mathbf{K}}_{0,t}^{(\varepsilon)}$, and since $\overline{J}_{0,t}^{(\varepsilon)}$ and $\overline{\mathbf{K}}_{0,t}^{(\varepsilon)}$ converge in ucp as $\varepsilon \to 0$ to some $\overline{J}_{0,t}$ and $\overline{\mathbf{K}}_{0,t}$ such that

$$\overline{\mathbf{Y}}_t = \overline{J}_{0,t} \overline{\mathbf{Y}}_0 + \overline{\mathbf{K}}_{0,t},$$

it follows that $\|\mathbf{Y}_t\|_{B,r} \leq \overline{\mathbf{Y}}_t$ for fixed $t \geq 0$, giving (9.1).

Similar quantities such as $\overline{J}_{s,t}^{(\varepsilon)}$ and $\overline{J}_{s,t}$ can be defined when going from time $s$ to time $t$, and similar results hold. Let $\overline{V}_t^{(\varepsilon)} := \alpha_0 \|\mathbf{e}\|_{B,r} + \|\mathbf{ea}'\|_{B,r} \overline{\mathbf{Y}}_{t-}^{(\varepsilon)}$ be the COGARCH(1,1) volatility corresponding to $\overline{\mathbf{Y}}^{(\varepsilon)}$. Define

$$X_t := -\lambda t - \sum_{0 < s \leq t} \log(1 + (\Delta L_s)^2 \|\mathbf{ea}'\|_{B,r}),$$

$$X_t^{(\varepsilon)} := -\lambda t - \sum_{0 < s \leq t, (\Delta L_s)^2 \geq \varepsilon} \log(1 + (\Delta L_s)^2 \|\mathbf{ea}'\|_{B,r}).$$

Then it follows from Theorem 2.2 and (1.4), that

$$\overline{V}_{t+}^{(\varepsilon)} = \left(\overline{V}_0 - \alpha_0 \|\mathbf{e}\|_{B,r} \lambda \int_0^t e^{X_s^{(\varepsilon)}} ds\right) e^{-X_t^{(\varepsilon)}}.$$

Thus we have $\overline{J}_{0,t}^{(\varepsilon)} = e^{-X_t^{(\varepsilon)}}$ and obtain another formula for $\overline{\mathbf{K}}_{0,t}^{(\varepsilon)}$, namely

$$\overline{\mathbf{K}}_{0,t}^{(\varepsilon)} = \|\mathbf{ea}'\|_{B,r}^{-1} \left[\alpha_0 \|\mathbf{e}\|_{B,r} e^{-X_t^{(\varepsilon)}} - \alpha_0 \|\mathbf{e}\|_{B,r} \lambda \int_0^t e^{-(X_t^{(\varepsilon)} - X_s^{(\varepsilon)})} ds - \alpha_0 \|\mathbf{e}\|_{B,r}\right].$$

From this it can be seen that $\overline{J}_{0,t}^{(\varepsilon)}$ and $\overline{\mathbf{K}}_{0,t}^{(\varepsilon)}$ are bounded by $\overline{J}_{0,t} = e^{-X_t}$ and

$$\overline{\mathbf{K}}_{0,t} = \|\mathbf{ea}'\|_{B,r}^{-1} \alpha_0 \|\mathbf{e}\|_{B,r} \left[e^{-X_t} - \lambda \int_0^t e^{-(X_t - X_s)} ds - 1\right],$$

respectively. Now define the versions

$$\overline{\mathbf{Y}}_\infty := \sum_{i=0}^\infty \overline{J}_{0,1} \cdots \overline{J}_{i-1,i} \overline{\mathbf{K}}_{i,i+1},$$

$$\mathbf{Y}_\infty^{(\varepsilon)} := \sum_{i=0}^\infty J_{0,1}^{(\varepsilon)} \cdots J_{i-1,i}^{(\varepsilon)} \mathbf{K}_{i,i+1}^{(\varepsilon)},$$

$$\mathbf{Y}_\infty := \sum_{i=0}^\infty J_{0,1} \cdots J_{i-1,i} \mathbf{K}_{i,i+1}.$$

In the proof of Theorems 3.1 and 3.3 we have seen that (3.2) implies that the sum defining $\overline{\mathbf{Y}}_\infty$ converges almost surely. This then gives the claim, since

$$\|J_{i-1,i}\|_{B,r}, \|J_{i-1,i}^{(\varepsilon)}\|_{B,r} \leq \overline{J}_{i-1,i},$$

$$\|\mathbf{K}_{i,i+1}\|_{B,r}, \|\mathbf{K}_{i,i+1}^{(\varepsilon)}\|_{B,r} \leq \overline{\mathbf{K}}_{i,i+1},$$

and $J_{i-1,i}^{(\varepsilon)}$ and $\mathbf{K}_{i,i+1}^{(\varepsilon)}$ converge in probability to $J_{i-1,i}$ and $\mathbf{K}_{i,i+1}$ as $\varepsilon \to 0$, respectively. $\square$

PROOF OF PROPOSITION 4.1. All assertions apart from the implication "(4.2) $\Longrightarrow \lambda(\widetilde{B}) < 0$" follow immediately from Lemma 9.1 (observing that the existence of $E\|Y_t\|^k$ is independent of the specific matrix norm) and the corresponding properties of the COGARCH$(1,1)$ process; see Section 4 in [19]. That (4.2) implies $\lambda(\widetilde{B}) < 0$ is a consequence of the Bauer–Fike perturbation result on eigenvalues, stating that for every eigenvalue $\widetilde{\lambda}_j$ of $\widetilde{B}$ we have

$$\min_{i=1,\ldots,q} |\lambda_i - \widetilde{\lambda}_j| \leq \|S^{-1}(\widetilde{B} - B)S\|_r = \mu \|S^{-1}\mathbf{e}\mathbf{a}'S\|_r$$

(see, e.g., Theorem 7.2.2 and its proof in [17]). $\square$

PROOF OF THEOREM 4.2. Since for fixed $t$, almost surely $V_t = V_{t+} = \alpha_0 + \mathbf{a}'\mathbf{Y}_t$, we obtain

(9.3) $$\operatorname{cov}(V_{t+h}, V_t) = \mathbf{a}' \operatorname{cov}(\mathbf{Y}_{t+h}, \mathbf{Y}_t) \mathbf{a}.$$

For ease of notation, we will assume that $t = 0$. Let $J_h := J_{0,h}$ and $\mathbf{K}_h := \mathbf{K}_{0,h}$ as constructed in the proof of Theorem 3.3. Then, using that $\|e^{Bt}\| \leq e^{\|B\|t}$ for any vector norm $\|\cdot\|$, it follows as in the proof of (8.6) that

(9.4) $$E\|J_h\| \leq e^{\|B\|t} E\left\{\exp\left(\sum_{0<s\leq h} \log(1+(\Delta L_s)^2\|\mathbf{e}\mathbf{a}'\|)\right)\right\} < \infty$$

by Klüppelberg, Lindner and Maller ([19], Lemma 4.1(a)). Using that $\mathbf{Y}_h = J_h \mathbf{Y}_0 + \mathbf{K}_h$, we conclude that $E\|\mathbf{K}_h\| < \infty$ and that

$$E(\mathbf{Y}_h \mathbf{Y}_0') = E(E(\mathbf{Y}_h \mathbf{Y}_0' | J_h, \mathbf{K}_h))$$
$$= E(J_h E(\mathbf{Y}_0 \mathbf{Y}_0') + \mathbf{K}_h E(\mathbf{Y}_0'))$$
$$= E(J_h) E(\mathbf{Y}_0 \mathbf{Y}_0') + E(\mathbf{K}_h) E(\mathbf{Y}_0').$$

On the other hand,
$$E(\mathbf{Y}_h) E(\mathbf{Y}_0') = E(J_h) E(\mathbf{Y}_0) E(\mathbf{Y}_0') + E(\mathbf{K}_h) E(\mathbf{Y}_0'),$$
so that $\text{cov}(\mathbf{Y}_h, \mathbf{Y}_0) = E(J_h) \text{cov}(\mathbf{Y}_0)$ and (4.4) will follow from (9.3) once we have shown that

(9.5) $$E(J_t) = e^{\widetilde{B}t}, \qquad t \geq 0.$$

To do that, it suffices to assume that $[L, L]_t$ is a compound-Poisson process. The general case then follows from the fact that $J_t^{(\varepsilon)}$ as defined in the proof of Theorem 3.1 converges to $J_t$ in $L^1$ as $\varepsilon \to 0$, since it converges stochastically and since there is an integrable majorant by (9.4) and its proof. So suppose that $[L, L]_t = \sum_{i=1}^{N(t)} Z_i$ is compound-Poisson with intensity $c > 0$ and let $C_i = (I + Z_i \mathbf{e}\mathbf{a}') e^{B(\Gamma_i - \Gamma_{i-1})}$. Then, for $0 \leq s, t$, it follows from (3.4) and the independence of $J_{0,s}$ and $J_{s,s+t}$ that
$$E(J_{s+t}) = E(J_s) E(J_t).$$
It is easy to see that $E(J_t)$ is a continuous function in $t \in [0, \infty)$. Furthermore, $E(J_0) = I$ and we conclude that $(E(J_t))_{t \geq 0}$ is a semigroup. We shall show that its generator $A_J$ satisfies

(9.6) $$A_J := \lim_{t \to 0} \frac{1}{t}(E(J_t) - I) = B + \int_{\mathbb{R}} y^2 \, d\nu_L(y) \, \mathbf{e}\mathbf{a}' = \widetilde{B}.$$

This then implies (9.5), since $E(J_t) = e^{tA_J}$ (see, e.g., [16], Proposition 2.5). To show (9.6), write

(9.7) $$\begin{aligned} J_t &= e^{Bt} \mathbf{1}_{\{N(t)=0\}} + e^{B(t-\Gamma_1)} C_1 \mathbf{1}_{\{N(t)=1\}} \\ &\quad + e^{B(t-\Gamma_{N(t)})} C_{N(t)} \cdots C_1 \mathbf{1}_{\{N(t) \geq 2\}}. \end{aligned}$$

Since $N(t)$ is Poisson distributed with parameter $ct$, we have $P(N(t) = k) = e^{-ct}(ct)^k/(k!)$. Then by (9.4),

$$E(e^{B(t-\Gamma_{N(t)})} C_{N(t)} \cdots C_1 \mathbf{1}_{\{N(t) \geq 2\}})$$

$$\leq e^{\|B\|t} E\left(\exp\left(\sum_{i=1}^{N(t)} \log(1 + Z_i \|\mathbf{e}\mathbf{a}'\|)\right) \mathbf{1}_{\{N(t) \geq 2\}}\right)$$

$$= e^{\|B\|t} E\left(\exp\left(\sum_{i=1}^{N(t)} \log(1 + Z_i \|\mathbf{e}\mathbf{a}'\|)\right) \Big| N(t) \geq 2\right) P(N(t) \geq 2)$$

(9.8) $$\leq e^{\|B\|t} E\left(\exp\left(\sum_{i=1}^{N(t)+2} \log(1 + Z_i \|\mathbf{e}\mathbf{a}'\|)\right)\right) P(N(t) \geq 2)$$

$$= e^{\|B\|t} E((1 + Z_1 \|\mathbf{e}\mathbf{a}'\|)(1 + Z_2 \|\mathbf{e}\mathbf{a}'\|))$$

$$\times E\left(\exp\left(\sum_{0<s\leq t}\log(1+(\Delta L_s)^2\|\mathbf{ea}'\|)\right)\right)P(N(t)\geq 2)$$

$$= o(t) \quad \text{as } t \to 0,$$

since $P(N(t) \geq 2) = o(t)$ as $t \to 0$. Furthermore, since $\Gamma_1$ is uniformly distributed on $(0, t)$, conditional on $N(t) = 1$, it follows that

$$E(e^{B(t-\Gamma_1)}C_1 \mathbf{1}_{\{N(t)=1\}})$$
$$= E(e^{B(t-\Gamma_1)}(I + Z_1\mathbf{ea}')e^{B\Gamma_1}|N(t)=1)P(N(t)=1)$$
$$= \int_0^t e^{B(t-s)}(I + E(Z_1)\mathbf{ea}')e^{Bs}\frac{ds}{t}e^{-ct}ct.$$

Since $\sup_{0 \leq s \leq t} \|e^{Bs} - I\|$ converges to 0 as $t \to 0$, we conclude that

$$\lim_{t \to 0} \frac{1}{t} E(e^{B(t-\Gamma_1)}C_1 \mathbf{1}_{\{N(t)=1\}}) = (I + E(Z_1)\mathbf{ea}')c.$$

Now (9.7) and (9.8) give (9.6), since

$$\lim_{t\to 0}\frac{E(J_t) - I}{t} = \lim_{t\to 0}\frac{e^{Bt}e^{-ct} - I}{t} + c(I + E(Z_1)\mathbf{ea}')$$
$$= -cI + B + c(I + E(Z_1)\mathbf{ea}') = \widetilde{B}. \qquad \square$$

We now need the following lemma:

LEMMA 9.2. *Let $T$ be exponentially distributed with parameter $c$, and suppose that $\lambda(B) < 0$. Let*

$$M := E(e^{BT} \otimes e^{BT}).$$

*Then*

(9.9) $\qquad E(e^{BT}) = (I - c^{-1}B)^{-1},$

(9.10) $\qquad M^{-1} = I_{q^2} - (I \otimes (c^{-1}B)) - ((c^{-1}B) \otimes I).$

*Furthermore, $(I \otimes B) + (B \otimes I)$ is invertible and, for any real $(q \times q)$ matrix $U$, the unique solution of $((I \otimes B) + (B \otimes I))\mathbf{x} = \text{vec}(U)$ is given by*

(9.11) $\qquad \mathbf{x} = \text{vec}\left(-\int_0^\infty e^{Bt}Ue^{B't}\,dt\right).$

*Here, we denote by $I$ the $(q \times q)$ identity matrix and denote by $I_{q^2}$ the $(q^2 \times q^2)$ identity matrix.*

PROOF. Equations (9.9) and (9.10) follow by simple calculations and a diagonalization argument, while invertibility of $(I \otimes B) + (B \otimes I)$ and (9.11) is a consequence of Lyapunov's theorem for the solution of Lyapunov equations (see, e.g., Section 9.3 in [15]). $\square$

PROOF OF LEMMA 4.3. Suppose first that the Lévy measure of $L$ is finite, and let $Q$ and $\mathbf{R}$ be as in Theorem 3.5(b) [writing $(T, Z)$ instead of $(T_0, Z_0)$]. Then, by Lemma 9.2,

$$E(Q) = (I - c^{-1}B)^{-1}(I + E(Z)\mathbf{ea}'),$$
$$E(\mathbf{R}) = \alpha_0 E(Z)(I - c^{-1}B)^{-1}\mathbf{e},$$

and (3.11) gives

$$(I - E(Q))E(\mathbf{Y}_\infty) = E(\mathbf{R}).$$

Furthermore,

$$(I - c^{-1}B)(I - E(Q)) = [(I - c^{-1}B) - I - E(Z)\mathbf{ea}']$$
$$= -\frac{1}{c}(B + \mu\mathbf{ea}'),$$

giving

$$E(\mathbf{Y}_\infty) = -c(B + \mu\mathbf{ea}')^{-1}(I - c^{-1}B)E(\mathbf{R}) = -\alpha_0\mu(B + \mu\mathbf{ea}')^{-1}\mathbf{e}.$$

Denoting $\mathbf{u} = (u_1, \ldots, u_q)' := (B + \mu\mathbf{ea}')^{-1}\mathbf{e}$, it is easy to see that $u_2 = \cdots = u_q = 0$ and $u_1 = 1/(\alpha_1\mu - \beta_q)$. In the case when $\nu_L$ is infinite, the result follows from Lemma 9.1, using that $\overline{Y}_\infty$ is an integrable majorant by (4.2). $\square$

PROOF OF THEOREM 4.4. By Lemma 9.1 and the dominated convergence theorem, for showing (4.6) it is sufficient to assume that $[L, L]$ is a compound-Poisson process. Hence, let $Q$ and $\mathbf{R}$ be as in Theorem 3.5, writing $(T, Z)$ instead of $(T_0, Z_0)$, where $T$ is exponentially distributed with parameter $c > 0$. Then

$$(9.12) \quad \begin{aligned} E(\mathbf{Y}_\infty \mathbf{Y}'_\infty) - E(Q\mathbf{Y}_\infty \mathbf{Y}'_\infty Q') \\ = E(Q\mathbf{Y}_\infty \mathbf{R}') + E(\mathbf{R}\mathbf{Y}'_\infty Q') + E(\mathbf{R}\mathbf{R}') \end{aligned}$$

by (3.11) and all these expectations exist by (4.3). Now

$$\begin{aligned} E(Q\mathbf{Y}_\infty \mathbf{Y}'_\infty Q') &= E(E[Q\mathbf{Y}_\infty \mathbf{Y}'_\infty Q'|Q]) \\ &= E(E[Q\,E(\mathbf{Y}_\infty \mathbf{Y}'_\infty)\,Q'|T]) \\ &= E(e^{BT} E[(I + Z\mathbf{ea}')E(\mathbf{Y}_\infty \mathbf{Y}'_\infty)(I + Z\mathbf{ae}')]e^{B'T}). \end{aligned}$$

Using that $\text{vec}(A_1 A_2 A_3) = (A_3' \otimes A_1)\text{vec}(A_2)$ for matrices $A_1, A_2$ and $A_3$, it follows with $M$ as in Lemma 9.2 that

$$\text{vec}(E(Q\mathbf{Y}_\infty \mathbf{Y}_\infty' Q'))$$
$$= M\text{vec}(E((I + Z\mathbf{ea}')E(\mathbf{Y}_\infty \mathbf{Y}_\infty')(I + Z\mathbf{ae}')))$$
$$= M(E((I + Z\mathbf{ea}') \otimes (I + Z\mathbf{ea}')))\text{vec}(E(\mathbf{Y}_\infty \mathbf{Y}_\infty'))$$
$$= M(I_{q^2} + E(Z)((\mathbf{ea}') \otimes I) + E(Z)(I \otimes (\mathbf{ea}')) + E(Z^2)((\mathbf{ea}') \otimes (\mathbf{ea}')))$$
$$\times \text{vec}(E(\mathbf{Y}_\infty \mathbf{Y}_\infty')).$$

Similar expressions can be obtained for $\text{vec}(E(Q\mathbf{Y}_\infty \mathbf{R}'))$, $\text{vec}(E(\mathbf{R}\mathbf{Y}_\infty' Q'))$ and $\text{vec}(E(\mathbf{R}\mathbf{R}'))$, and we obtain from (9.12) that

$$[I_{q^2} - M(I_{q^2} + E(Z)((\mathbf{ea}') \otimes I) + E(Z)(I \otimes (\mathbf{ea}')) + E(Z^2)((\mathbf{ea}') \otimes (\mathbf{ea}')))]$$
$$\times \text{vec}(E(\mathbf{Y}_\infty \mathbf{Y}_\infty'))$$
$$= M\text{vec}[\alpha_0^2 E(Z^2)\mathbf{ee}' + \alpha_0(E(Z)I + E(Z^2)\mathbf{ea}')E(\mathbf{Y}_\infty)\mathbf{e}'$$
$$+ \alpha_0 \mathbf{e}E(\mathbf{Y}_\infty')(E(Z)I + E(Z^2)\mathbf{ae}')].$$

Multiplying this equation by $cM^{-1}$, using (9.10) and (4.5) as well as $\mu = cE(Z)$ and $\rho = cE(Z^2)$, we obtain

$$-[(I \otimes (B + \mu\mathbf{ea}')) + ((B + \mu\mathbf{ea}') \otimes I) + \rho((\mathbf{ea}') \otimes (\mathbf{ea}'))]\text{vec}(E(\mathbf{Y}_\infty \mathbf{Y}_\infty'))$$
$$= \text{vec}[\alpha_0^2 \rho\mathbf{ee}' - \alpha_0^2(\mu I + \rho\mathbf{ea}')\mu(B + \mu\mathbf{ea}')^{-1}\mathbf{ee}'$$
$$- \alpha_0^2 \mathbf{ee}'(B' + \mu\mathbf{ae}')^{-1}\mu(\mu I + \rho\mathbf{ae}')].$$

Adding to this

$$[(I \otimes \widetilde{B}) + (\widetilde{B} \otimes I) + \rho((\mathbf{ea}') \otimes (\mathbf{ea}'))]\text{vec}(E(\mathbf{Y}_\infty)E(\mathbf{Y}_\infty'))$$
$$= \text{vec}[\widetilde{B}\,E(\mathbf{Y}_\infty)E(\mathbf{Y}_\infty') + E(\mathbf{Y}_\infty)E(\mathbf{Y}_\infty')\widetilde{B}' + \rho\mathbf{ea}'E(\mathbf{Y}_\infty)E(\mathbf{Y}_\infty')\mathbf{ae}']$$
$$= \alpha_0^2 \text{vec}[\mu^2 \mathbf{ee}'(\widetilde{B}')^{-1} + \mu^2 \widetilde{B}^{-1}\mathbf{ee}' + \rho\mu^2 \mathbf{ea}'\widetilde{B}^{-1}\mathbf{ee}'(\widetilde{B}')^{-1}\mathbf{ae}']$$

on both sides results in

$$-[(I \otimes \widetilde{B}) + (\widetilde{B} \otimes I) + \rho((\mathbf{ea}') \otimes (\mathbf{ea}'))]\text{vec}(\text{cov}(\mathbf{Y}_0))$$
$$= \alpha_0^2 \rho[1 - \mu(\mathbf{a}'\widetilde{B}^{-1}\mathbf{e})]^2 \text{vec}(\mathbf{ee}')$$
$$= \frac{\alpha_0^2 \beta_q^2 \rho}{(\beta_q - \mu\alpha_1)^2}\text{vec}(\mathbf{ee}'),$$

which is (4.6), where we used (4.5) in the last equation.

Now let $A := (I \otimes \widetilde{B}) + (\widetilde{B} \otimes I)$ and $\mathbf{x} := \text{vec}(\text{cov}(\mathbf{Y}_\infty))$. By Proposition 4.1 and Lemma 9.2, $A$ is invertible. Observe that the matrix $\rho((\mathbf{ea}') \otimes (\mathbf{ea}'))$ has nonzero entries only in the last row. Denote this row by $\mathbf{c}'$. Furthermore, set $\gamma := \rho \alpha_0^2 \beta_q^2 (\mu \alpha_1 - \beta_q)^{-2}$. Then (4.6) can be written as

$$A\mathbf{x} + (\mathbf{c}'\mathbf{x})\mathbf{e}_{q^2} = -\gamma \, \mathbf{e}_{q^2}.$$

We know already that a solution to this equation exists. Suppose there are two solutions, call them $\mathbf{x}_1$ and $\mathbf{x}_2$. Then $A\mathbf{x}_1 = -(\gamma + \mathbf{c}'\mathbf{x}_1)\mathbf{e}_{q^2}$ and $A\mathbf{x}_2 = -(\gamma + \mathbf{c}'\mathbf{x}_2)\mathbf{e}_{q^2}$. Denoting the unique solution of $A\mathbf{y} = -n\,\mathbf{e}_{q^2}$ by $\mathbf{y}(n)$, $n \in \mathbb{R}$, it follows that $\mathbf{x}_1 = \mathbf{y}(\gamma + \mathbf{c}'\mathbf{x}_1)$ and $\mathbf{x}_2 = \mathbf{y}(\gamma + \mathbf{c}'\mathbf{x}_2)$. Since $\mathbf{x}_1 \neq \mathbf{0} \neq \mathbf{x}_2$, this implies $\gamma + \mathbf{c}'\mathbf{x}_1 \neq \mathbf{0} \neq \gamma + \mathbf{c}'\mathbf{x}_2$, and using the linearity of the solution $\mathbf{y}(n)$ in $n$, it follows that there is $\kappa \neq 0$ such that $\mathbf{x}_2 = \kappa \mathbf{x}_1$. Thus we have $A\mathbf{x}_1 = -(\gamma + \mathbf{c}'\mathbf{x}_1)\mathbf{e}_{q^2}$ and $\kappa A\mathbf{x}_1 = -(\gamma + \kappa \mathbf{c}'\mathbf{x}_1)\mathbf{e}_{q^2}$, and this is only possible if $\kappa = 1$, so $\mathbf{x}_1 = \mathbf{x}_2$. So the solution of (4.6) is unique, implying that the matrix $A + \rho((\mathbf{e}\mathbf{a}') \otimes (\mathbf{e}\mathbf{a}'))$ is invertible.

By (9.11), the solution $y(n)$ of $Ay = -ne_{q^2}$ is given by

$$(9.13) \qquad y(n) = \mathrm{vec}\left( n \int_0^\infty e^{\widetilde{B}t} \mathbf{e}\mathbf{e}' e^{\widetilde{B}'t} \, dt \right).$$

This gives

$$\mathrm{cov}(\mathbf{Y}_\infty) = (\gamma + \mathbf{c}' \, \mathrm{vec}(\mathrm{cov}(\mathbf{Y}_\infty))) \int_0^\infty e^{\widetilde{B}t} \mathbf{e}\mathbf{e}' e^{\widetilde{B}'t} \, dt.$$

Since both $\mathrm{cov}(\mathbf{Y}_\infty)$ and $\int_0^\infty e^{\widetilde{B}t} \mathbf{e}\mathbf{e}' e^{\widetilde{B}'t} \, dt$ are positive semidefinite, it follows that $\gamma + \mathbf{c}' \mathrm{vec}(\mathrm{cov}(\mathbf{Y}_\infty)) > 0$. By Brockwell [9], the stationary CARMA state vector $\zeta_\infty$ has covariance matrix

$$\mathrm{cov}(\zeta_\infty) = \rho \int_0^\infty e^{\widetilde{B}t} \mathbf{e}\mathbf{e}' e^{\widetilde{B}'t} \, dt,$$

so that there is $u > 0$ such that

$$(9.14) \qquad \mathrm{cov}(\mathbf{Y}_\infty) = u \, \mathrm{cov}(\zeta_\infty).$$

Inserting (9.14) into (4.6) and using (9.13) shows

$$-u\rho \, \mathrm{vec}(\mathbf{e}\mathbf{e}') + u\rho^2 \, \mathrm{vec}\left( \mathbf{e}\mathbf{a}' \int_0^\infty e^{\widetilde{B}t} \mathbf{e}\mathbf{e}' e^{\widetilde{B}'t} \, dt \, \mathbf{a}\mathbf{e}' \right)$$
$$= \frac{-\alpha_0^2 \beta_q^2 \rho}{(\beta_q - \mu\alpha_1)^2} \, \mathrm{vec}(\mathbf{e}\mathbf{e}'),$$

so that

$$-u(1-m) \, \mathrm{vec}(\mathbf{e}\mathbf{e}') = \frac{-\alpha_0^2 \beta_q^2}{(\beta_q - \mu\alpha_1)^2} \, \mathrm{vec}(\mathbf{e}\mathbf{e}').$$

Since $u > 0$ and $\alpha_0, \beta_q \neq 0$, it follows that $0 \leq m < 1$ and that

$$u = \frac{\alpha_0^2 \beta_q^2}{(\beta_q - \mu\alpha_1)^2 (1-m)},$$

giving (4.7). This implies (4.8), using $V_\infty = \alpha_0 + \mathbf{a}'\mathbf{Y}_\infty$, and (4.9) follows from (4.5). Finally,

$$E(\psi_\infty) = \mathbf{a}'E\int_0^\infty e^{\widetilde{B}t}\mathbf{e}\,d\widetilde{L}_t = \mu\int_0^\infty \mathbf{a}'e^{\widetilde{B}t}\mathbf{e}\,dt = -\mu\mathbf{a}'\widetilde{B}^{-1}\mathbf{e},$$

giving (4.10), and (4.11) and (4.12) are direct consequences of (4.4), (4.7) and the autocovariance function of a CARMA process (see [9]). $\square$

## 10. Proofs for Section 5.

PROOF OF THEOREM 5.1. (a) Suppose that (5.1) and (5.2) both hold. By Lemma 8.2, it suffices to show (5.3) for the case that $[L,L] = \sum_{i=1}^{N(t)} Z_i$ is a compound-Poisson process, with jump times $(\Gamma_n)_{n\in\mathbb{N}}$. Then it follows easily by induction from (2.3) and (8.1) that

$$\mathbf{Y}_t = e^{Bt}\mathbf{Y}_0 + \sum_{i=1}^{N(t)} e^{B(t-\Gamma_i)}\mathbf{e}V_{\Gamma_i}Z_i, \qquad t \geq 0.$$

In view of the proof of (b) below, let $s \geq 0$. Then

$$(10.1) \qquad \mathbf{a}'e^{Bs}\mathbf{Y}_t = \mathbf{a}'e^{B(s+t)}\mathbf{Y}_0 + \sum_{i=1}^{N(t)} \mathbf{a}'e^{B(s+t-\Gamma_i)}\mathbf{e}V_{\Gamma_i}Z_i$$

$$(10.2) \qquad \geq \gamma + \sum_{i=1}^{N(t)} \mathbf{a}'e^{B(s+t-\Gamma_i)}\mathbf{e}V_{\Gamma_i}Z_i.$$

Setting $s = 0$, it follows that $V_t = \alpha_0 + \mathbf{a}'\mathbf{Y}_{t-} \geq \alpha_0 + \gamma$ for $t \in [0,\Gamma_1]$; hence also $V_{\Gamma_1+} \geq \alpha_0 + \gamma \geq 0$ by (5.1) and (10.2), and an induction argument shows that $V_t \geq \alpha_0 + \gamma$ for all $t \geq 0$, that is, (5.3) holds.

For the converse, suppose first that (5.2) fails. Then, using the continuity of the function $t \mapsto e^{Bt}$, it follows that there is $(t_1,t_2) \subset (0,\infty)$ such that $P(\alpha_0 + \mathbf{a}'e^{Bt}\mathbf{Y}_0 < 0 \;\forall t \in (t_1,t_2)) > 0$, and since $P(\Gamma_1 > t_2) > 0$ we get the claim from (10.1). So suppose that (5.2) holds with $\gamma > -\alpha_0$, but (5.1) fails. Suppose that the support of the Lévy measure of the compound-Poisson process $[L,L]$ (and hence the support distribution of the jumps, $Z_i$) is unbounded. Let $(t_3,t_4) \subset (0,\infty)$ be an interval such that $\mathbf{a}'e^{Bt}\mathbf{e} \leq -c_1 < 0$ for all $t \in (t_3,t_4)$ for some $c_1 < 0$. Let $t_5 > t_4$. By (5.2) we have $P(V_{\Gamma_1} \geq \alpha_0 + \gamma) = 1$, so that it is easy to see that the set

$$A := \{\Gamma_1 < t_5 < \Gamma_2, t_5 - \Gamma_1 \in (t_3,t_4), V_{\Gamma_1} \geq \alpha_0 + \gamma\}$$

has positive probability. On $A$ we have, by (10.1),

$$V_{t_5} = \alpha_0 + \mathbf{a}'e^{Bt_5}\mathbf{Y}_0 + \mathbf{a}'e^{B(t_5-\Gamma_1)}\mathbf{e}V_{\Gamma_1}Z_1.$$

Now $\mathbf{a}'e^{B(t_5-\Gamma_1)}\mathbf{e} \leq -c_1$ and by choosing $Z_1$ (which is independent of $\Gamma_1, \Gamma_2$ and $\mathbf{Y}_0$) large enough, we obtain $P(V_{t_5} < 0) > 0$.

(b) In view of (a) it remains to show that $\mathbf{Y}_\infty$ satisfies (5.2). For the proof of this, it suffices by Lemma 9.1 to assume that $[L, L]$ is compound-Poisson. Let $(\widetilde{\mathbf{Y}}_t)_{t\geq 0}$ be a state process with $\widetilde{\mathbf{Y}}_0 = 0$. Then (5.2) holds for $\widetilde{\mathbf{Y}}_0$ with $\gamma = 0$, and it follows from (10.2), (5.1) and (5.3) that $\mathbf{a}'e^{Bs}\widetilde{\mathbf{Y}}_t \geq 0$ for all $s, t \geq 0$. Since $\widetilde{\mathbf{Y}}_t$ converges in distribution to $\mathbf{Y}_\infty$ as $t \to \infty$, (5.2) follows with $\gamma = 0$. □

## 11. Proof for Section 6.

PROOF OF THEOREM 6.1. We mimic the proof of Proposition 5.1 of [19], that is, in the COGARCH(1,1) case. Observe that (6.1) and (6.3) follow immediately, since $(L_t)_{t\geq 0}$ is a zero-mean martingale. Furthermore, $(G_t)_{t\geq 0}$ is a square integrable martingale so that

$$EG_r^2 = E\int_0^r V_s\, d[L,L]_s = E(L_1)^2 r E(V_\infty),$$

and (6.2) follows from (4.9). Before showing (6.4), we verify that $EG_t^4 < \infty$ if (4.3) is satisfied: it follows from the Burkholder–Davis–Gundy inequality (see, e.g., [24], page 222) that $EG_t^4 < \infty$ if $E[G,G]_t^2 < \infty$. Let $\overline{V}_t = \alpha_0\|\mathbf{e}\|_{B,r} + \|\mathbf{ea}'\|_{B,r}\overline{\mathbf{Y}}_{t-}$ be the volatility of the COGARCH(1,1) process constructed in Lemma 9.1 and let $\overline{G}_t = \int_0^t \sqrt{\overline{V}_t}\, dL_t$ be the corresponding COGARCH(1,1) price process. Then it follows from Lemma 9.1 that there is $C_1 > 0$ such that

$$0 \leq V_s = \alpha_0 + \mathbf{a}'\mathbf{Y}_{s-} \leq \alpha_0 + C_1\overline{\mathbf{Y}}_{s-} = \alpha_0 + C_1\frac{\overline{V}_s - \alpha_0\|\mathbf{e}\|_{B,r}}{\|\mathbf{ea}'\|_{B,r}}.$$

Then

$$[G,G]_t = \int_0^t V_s\, d[L,L]_s$$
$$\leq \frac{C_1}{\|\mathbf{ea}'\|_{B,r}}\int_0^t \overline{V}_s\, d[L,L]_s + \left(\alpha_0 - \frac{C_1\alpha_0\|\mathbf{e}\|_{B,r}}{\|\mathbf{ea}'\|_{B,r}}\right)[L,L]_t$$
$$= \frac{C_1}{\|\mathbf{ea}'\|_{B,r}}[\overline{G},\overline{G}]_t + \left(\alpha_0 - \frac{C_1\alpha_0\|\mathbf{e}\|_{B,r}}{\|\mathbf{ea}'\|_{B,r}}\right)[L,L]_t,$$

so that again by the Burkholder–Davis–Gundy inequality and Doob's maximal inequality, finiteness of $E\overline{G}_t^4$ implies finiteness of $E[\overline{G},\overline{G}]_t^2$ and hence of $EG_t^4$. The fact that $E\overline{G}_t^4 < \infty$ was used in [19] in the case when $L$ has no Gaussian component, but it also holds in the general case.

Denote by $E_r$ the conditional expectation with respect to the $\sigma$-algebra $\mathcal{F}_r$. Using partial integration, we have

$$(G_h^{(r)})^2 = 2\int_{h+}^{h+r} G_{s-}\, dG_s + [G,G]_{h+}^{h+r}$$

$$= 2\int_h^{h+r} G_{s-}\sqrt{V_s}\, dL_s + \int_{h+}^{h+r} V_s\, d[L,L]_s.$$

Since the increments of $L$ on the interval $(h, h+r]$ are independent of $\mathcal{F}_r$ and since $L$ has expectation 0, it follows that

$$E_r \int_{h+}^{h+r} G_{s-}\sqrt{V_s}\, dL_s = 0.$$

Recall that $\mathbf{Y}_s = J_{r,s}\mathbf{Y}_r + \mathbf{K}_{r,s}$ by (3.3). Hence we also have $\mathbf{Y}_{s-} = J_{r,s-}\mathbf{Y}_r + \mathbf{K}_{r,s-}$, so that, by the compensation formula,

$$E_r(G_h^{(r)})^2 = E_r \int_{h+}^{h+r} (\alpha_0 + \mathbf{a}'\mathbf{Y}_{s-})\, d[L,L]$$

$$= E_r \int_{h+}^{h+r} (\alpha_0 + \mathbf{a}'J_{r,s-}\mathbf{Y}_r + \mathbf{a}'\mathbf{K}_{r,s-})\, d[L,L]$$

(11.1)
$$= E(L_1^2)\alpha_0 r + E(L_1^2)\mathbf{a}' \int_{h+}^{h+r} (EJ_{r,s-})\mathbf{Y}_r\, ds$$

$$\quad + E(L_1^2)\mathbf{a}' \int_{h+}^{h+r} (E\mathbf{K}_{r,s-})\, ds$$

$$= E(L_1^2) \int_h^{h+r} E_r(V_s)\, ds.$$

Since $\mathbf{Y}_\infty \stackrel{d}{=} J_{r,s}\mathbf{Y}_\infty + \mathbf{K}_{r,s}$ by (3.5), with $\mathbf{Y}_\infty$ independent of $(J_{r,s}, \mathbf{K}_{r,s})$ on the right-hand side, and $EJ_{r,s} = e^{\widetilde{B}(s-r)}$ by the proof of Theorem 4.2, it follows from (4.5) that

$$E\mathbf{K}_{r,s} = (I - e^{\widetilde{B}(s-r)})\frac{\alpha_0 \mu}{\beta_q - \alpha_1 \mu}\mathbf{e}_1.$$

Hence

(11.2)
$$E_r(V_s) = \alpha_0 + \mathbf{a}'e^{\widetilde{B}(s-r)}\mathbf{Y}_r + \mathbf{a}'\frac{\alpha_0 \mu}{\beta_q - \alpha_1 \mu}(I - e^{\widetilde{B}(r-s)})\mathbf{e}_1$$

$$= \frac{\alpha_0 \beta_q}{\beta_q - \alpha_1 \mu} + \mathbf{a}'e^{\widetilde{B}(s-r)}\left(\mathbf{Y}_r - \frac{\alpha_0 \mu}{\beta_q - \alpha_1 \mu}\mathbf{e}_1\right).$$

Combining $\int_h^{h+r} e^{\widetilde{B}(s-r)}\, ds = e^{\widetilde{B}h}\widetilde{B}^{-1}(I - e^{-\widetilde{B}r})$ with (11.1), (11.2) and (4.5) gives

$$E_r(G_h^{(r)})^2 = E(L_1^2)\left(\frac{\alpha_0 r \beta_q}{\beta_q - \alpha_1 \mu} + \mathbf{a}'e^{\widetilde{B}h}\widetilde{B}^{-1}(I - e^{-\widetilde{B}r})(\mathbf{Y}_r - E\mathbf{Y}_r)\right),$$

and we conclude with (6.2) that

$$\begin{aligned}
&E((G_0^{(r)})^2 (G_h^{(r)})^2) \\
&= E(E_r((G_h^{(r)})^2 G_r^2)) \\
&= E(L_1^2) E\left(\frac{\alpha_0 r \beta_q}{\beta_q - \alpha_1 \mu} G_r^2 + \mathbf{a}' e^{\widetilde{B}h} \widetilde{B}^{-1}(I - e^{-\widetilde{B}r})(\mathbf{Y}_r - E\mathbf{Y}_r) G_r^2\right) \\
&= (E(G_r^2))^2 + E(L_1^2) \mathbf{a}' e^{\widetilde{B}h} \widetilde{B}^{-1}(I - e^{-\widetilde{B}r})[E(\mathbf{Y}_r G_r^2) - (E\mathbf{Y}_r)(E G_r^2)],
\end{aligned}$$

showing (6.4). □

**Acknowledgments.** We are indebted to two referees for a number of helpful comments and corrections, and to Ole Barndorff-Nielsen, Vicky Fasen, Claudia Klüppelberg and Robert Stelzer for fruitful and inspiring comments. A. Lindner also thanks the Statistics Department at Colorado State University and P. Brockwell thanks the Mathematics Department at the Technical University of Munich and the Department of Mathematics and Statistics at the University of Melbourne for their stimulating and enjoyable hospitality.

P. Brockwell  
E. Chadraa  
Department of Statistics  
Colorado State University  
Fort Collins, Colorado 80523-1877  
USA  
E-mail: pjbrock@stat.colostate.edu  
    chadraa@stat.colostate.edu

A. Lindner  
Centre for Mathematical Sciences  
Technical University of Munich  
D-85747 Garching  
Germany  
E-mail: lindner@ma.tum.de